\documentclass[reqno]{amsart} 
\usepackage{amsmath,amssymb, xypic}

\numberwithin{equation}{section}

\newcommand\pone{\Bbb{P}^1}
\newcommand\mg{\mathcal{G}}
\newcommand\mf{\mathcal{F}}

\newcommand{\tw}{\operatorname{\tilde{w}}}
\newcommand{\tmu}{\operatorname{\tilde{\mu}}}
\newcommand{\tmw}{\operatorname{\tilde{\mathcal{W}}}}

\newcommand\la{{\lambda}}
\newcommand\mv{\mathcal{V}}
\newcommand\ma{\mathcal{A}}
\newcommand\ms{\mathcal{S}}
\newcommand\mt{\mathcal{T}}
\newcommand\mw{\mathcal{W}}
\newcommand\ml{\mathcal{L}}
\newcommand\mq{\mathcal{Q}}
\newcommand\muu{\mathcal{U}}
\newcommand\bull{\sssize{\bullet}}

\newtheorem{theorem}{Theorem}[section]
\newtheorem{lemma}[theorem]{Lemma}
\newtheorem{question}[theorem]{Question}
\newtheorem{proposition}[theorem]{Proposition}
\newtheorem{corollary}[theorem]{Corollary}
\newtheorem{definition-proposition}[theorem]{Definition-Proposition}

\theoremstyle{definition} \newtheorem{definition}[theorem]{Definition}
\theoremstyle{defi} \newtheorem{defi}[theorem]{Definition}
\newtheorem{example}[theorem]{Example}

 \newtheorem{claim}[theorem]{Claim}

\theoremstyle{remark} \newtheorem{remark}[theorem]{Remark}

\begin{document}
\title[Eigenvalue problems]{Irredundancy in Eigenvalue problems}
\author{Prakash Belkale}
\address{Department of Mathematics\\ UNC-Chapel Hill\\ CB \#3250, Phillips Hall
\\ Chapel Hill , NC 27599}
\email{belkale@email.unc.edu} 

\begin{abstract}
 We show the irredundancy, necessity and sufficiency of inequalities corresponding to  intersections  whose ``Witten Bundle'' is polyrigid for the product of unitary matrices problem  ~\cite{AW},~\cite{b1}.  The irredundancy problem in the classical case (in a stronger form) is due to Knutson, Tao and Woodward. Our methods (specialized to the classical case) are different from ~\cite{KTW}.

\end{abstract}
\maketitle
\section*[intro]{Introduction}

In ~\cite{AW} and ~\cite{b1} the following problem was considered. Let
$\bar{A}_1,\dots,\bar{A}_s$, be conjugacy classes in $SU(n)$, find conditions on these so that there are elements $\bar{A}_i\in SU(n)$ with conjugacy class of
$A_i$ in $\bar{A}_i$ and

$$A_1A_2\dots A_s=1$$

This problem is  in a natural way related to unitary representations
of the fundamental group of $\Bbb{P}^1-\{p_1,\dots,p_s\}$ where $\{p_1,\dots,p_s\}$ is a set of distinct points on $\Bbb{P}^1$. The answer to this question involves Quantum Cohomology of Grassmannians.  Recall the principal result in ~\cite{b1}:

\begin{theorem}Let $\bar{A}_l,l=1,\dots,s$ be conjugacy classes in $SU(n)$.
. Then there exist $A_l\in SU(n)$with conjugacy class of $A_i$ in $\bar{A}_i$ and $A_1A_2\dots A_s=1$ if and only if  
given any integers $1\leq r<n, d\geq 0$,subsets $I^1,\dots I^s$ of $\{1,\dots,n\}$ each of cardinality $r$, such that the Gromov-Witten number
$$\langle\sigma_{I^1},\dots,\sigma_{I^{s}}\rangle_{d}=1,$$

the inequality\footnote{See section ~\ref{Evalue1} for the notation, and the definition of $\lambda_{I}(\bar{A})$.}

$$\sum_{l=1}^{l=s}\lambda_{I^l}(\bar{A}_l) -d\leq 0$$

holds. 
\end{theorem}

Let us name the inequality corresponding to the data  $(d,r,I^1,\dots,I^s)$
as\linebreak $Ineq(d,r,I^1,\dots,I^s)$. So according to the theorem above the set of inequalities $Ineq(d,r,I^1,\dots,I^s)$ satisfying
$$\langle\sigma_{I^1},\dots,\sigma_{I^{s}}\rangle_{d}=1$$

form a necessary and sufficient set of inequalities defining the desired
set.

The reduction to the set of inequalities corresponding to intersection $=1$
was made using the existence of Harder-Narasimhan filtration. The natural question at this point is whether there are more such conditions and to find an irredundant (and sufficient) set of inequalities for this problem. Showing irredundancy of an inequality entails the  construction (and is equivalent to it\footnote{Lemma ~\ref{lemmacondition}.} ) of regular\footnote{That is, having distinct eigenvalues.} elements $A_i$ in SU(n) with product $=I$ for which a given inequality turns into an equality. So we need a source of unitary representations of $\pi_1(\Bbb{P}^1-\{p_1,\dots,p_s\})$ to prove irredundancy statements.

Recall the theorem of Witten which associates a moduli space of parabolic bundles on $\Bbb{P}^1$
to a GW number (these results are summarized in ~\ref{HWitten}). 

We describe the main results of this paper:
\begin{enumerate}
\item There are more such conditions. We show that if $Ineq(d,r,I^1,\dots,I^s)$
 is a necessary inequality \footnote{ In fact if the subbundle  corresponding to it is the Harder-Narasimhan element for some choice of weights} the  then the associated Witten moduli space of parabolic bundles is a point. 
\item We show that the $Ineq(d,r,I^1,\dots,I^s)$ which have associated Witten moduli space a point are irredundant. In fact we produce weights for which the subbundle  corresponding to $Ineq(d,r,I^1,\dots,I^s)$ is the only contradiction to semistability (among all subbundles).
\end{enumerate} 
\begin{defi}
We will call an inequality $Ineq(d,r,I^1,\dots,I^s)$ for which the associated Witten moduli space is a point to be a polyrigid inequality. That is, the ``associated'' unitary representation of the fundamental group of $\Bbb{P}^1-\{p_1,\dots,p_s\}$ is a direct sum of rigid representations (in the sense of N. Katz ~\cite{KATZ}). This is a computable property (see ~\ref{decision}).
\end{defi}

\begin{remark}
It is not known whether 
$$\langle\sigma_{I^1},\dots,\sigma_{I^{s}}\rangle_{d}=1,$$

implies that the associated Witten moduli space is a point. We will call this the Intersection rigidity question (IRQ). The classical part of this is Fulton's conjecture (theorem of Knutson-Tao-Woodward). 
\end{remark}

In the basic form the idea for showing irredundancy is simple: Suppose
$$\langle\sigma_{I^1},\dots,\sigma_{I^{s}}\rangle_{d}=1$$ and is polyrigid.
We get a parabolic structure on the point of intersection $\mv\subset \mt=\mathcal{O}^{\oplus n}$ by the easy part of the Quantum Horn conjecture. We can extend this structure to $\mt$ in a minimal way (see definition ~\ref{min} for the definition of minimal extension) and hope that this should be a structure which is not semistable and $\mv$ is the only contradiction to semistability. This works out in this form if both $\mv$ and its Grassmann dual are stable. For the rest of the cases we need an induction apparatus.

In section ~\ref{GW} we recall the setting for IRQ. We recall the  definition (which has appeared in ~\cite{b4}) of  generalized Gromov-Witten numbers of Grassmannians
(instead of looking for subbundles of the trivial bundle we generalize to subbundles of ``Generic bundles''). Computationally these can be recovered from the usual GW numbers.

In section ~\ref{GDuality} for future reference we make explicit Grassmann Duality for these generalized Gromov-Witten numbers.

In section ~\ref{PrePar} we make the most ``liberated'' definition of parabolic bundles with as few conditions as possible (the weights need not be between $0$ and $1$, but just that the pairwise difference of weights at any point is $\leq 1$ in absolute value). And we recall the construction of  Harder and Narasimhan
~\cite{HN} of maximally destabilising subbundles with this notation.

 In section ~\ref{HWitten} we recall our construction from ~\cite{b4} of
a pre-parabolic bundle from points of intersection corresponding to a generalized Gromov-Witten number. It is the easy case of the Quantum Horn conjecture that this bundle is semistable. We also recall the theorem of Witten which shows that GW numbers are the dimensions of global sections of certain line bundles over
suitable moduli spaces (this is the quantum generalization of the connection between representation theory of $GL(r)$ and intersection theory on Grassmannnians). In section 1.5 we have a warm up for the geometric study of  saturation phenomena, by a geometric proof of the easy part of the (quantum) saturation conjecture. This setting is needed for the proof of redundancy of inequalities which
violate (IRQ). 

In section ~\ref{Evalue1} we recall the main theorems from ~\cite{AW},~\cite{b1} on eigenvalue problems and state them in a form applicable to pre-parabolic bundles. We make the important definition of `polyrigid inequalities'. These are inequalities coming from intersections whose associated Witten \footnote{ The easy part of Quantum Horn conjecture essentially follows immediately from a theorem of Witten, but the actual construction of a bundle corresponding to a quantum intersection is in ~\cite{b4}.} moduli space is a point. We first discuss what redundancy means, and prove a condition equivalent to the irredundancy of an inequality.

In section ~\ref{redundancy} we show that the inequalities that are not polyrigid are redundant. In fact there do not exist any choice of weights for which they are a part of the Harder-Narasimhan Filtration (not just the first one). 

In Section ~\ref{irred} we compute dimensions of parabolic moduli spaces
under the assumption that there is a stable point in the moduli space. Now assuming that the moduli space is a point\footnote{ That is assuming rigidity.} (it is always connected and irreducible)
we construct another remarkable choice of weights (called rigidity weights) on the same parabolic bundle. We show that for this new weights the structure is
stable. 

The true origin of this structure is the following: Assuming
$$\langle\sigma_{I^1},\dots,\sigma_{I^{s}}\rangle_{d}\neq 0$$ and the Grassmann
duals of these cycles normalised (see definition ~\ref{NORMAL}), and with notation as above, we have a parabolic structure on $\mv$, which we extend to $\mt$
in the minimal possible way. hence there is a parabolic structure on $\mq=\frac{\mt}{\mv}$. These turn out to be the  rigidity weights for the Grassmann dual parabolic structure (which is on $\mq^*$, and hence on $\mq$).

In section ~\ref{I1} we start working towards results on Irredundancy. We show that irredundancy of the polyrigid inequalities follows from a structural theorem on polyrigid intersections (theorem ~\ref{UneTheoreme}).

In section ~\ref{examples} we give examples for the irredundancy construction of section ~\ref{I1}.
 
The proof of theorem ~\ref{UneTheoreme} is carried out in section ~\ref{JHSaga}. We establish a general theorem on Jordan Holder filtrations and use this to get an inductive grip on polyrigid inequalities,

Our proof of irredundancy of the polyrigid inequalities is by (almost explicitly) displaying situations were all but one inequality is true. And these situations come from the easy part of the Quantum Horn conjecture. We are unable to compare our `examples' of irredundancy with those in ~\cite{KTW} (or in ~\cite{DW2}). 

Although IRQ remains an important open problem, the problem of a provably irredundant (and sufficient) set of inequalities for the eigenvalues of product of unitary matrices has been solved. IRQ is a problem on the intersection theory side.

We would like to acknowledge the influence\footnote{We have not used any of their results, but the possibility of obtaining irrdundancy results which is a priori a bit far fetched comes from their work.}  of ~\cite{KTW},~\cite{DW2} on this work. The classical part of redundancy of inequalities which violate IRQ appears in ~\cite{DW2} in the language of quivers. They prove this statement\footnote{ The classical part of IRQ is Fulton's conjecture, which is a theorem of Knutson, Tao and Woodward.} by methods which seem different from ours. We  thank S. Kumar, C. Woodward for some useful conversations and encouragement. The author was partially supported by NSF grant DMS-0300356.

\section{Some generalities}

\subsection{Generalized Gromov-Witten numbers}\label{GW}
We recall the definition of Gromov-Witten invariants and also define some more general numbers.

      Let $I^{1},I^{2},\dots,I^{s}$ be subsets of $\{1,2,\dots,n\}$ of cardinality $r$ and $p_{1},\dots,p_{s}$ general points of $ {\mathbb P}^{1}$. Let
$\mathcal{W}=\mathcal{O}_{\Bbb{P}^1}^n$. The Gromov-Witten number

$$\langle\sigma_{I^1},\dots,\sigma_{I^{s}}\rangle_{d}$$

 is defined to be, for generic flags $F(\mw_{p_l})_{\bull}$ on $\mw_{p_l}$ ($l = 1,2,\dots,s$), the number of subbundles $\mathcal{V}$ of $\mw$ whose
fiber $\mathcal{V}_{p_l}\subset \mw_{p_l}$ lies in the Schubert cell $\Omega^{o}_{I^l}(F(\mw_{p_l})_{\bull}), l=1,\dots,s $ (the number is defined to be zero if there is a positive dimensional family of these).

Generic bundles on $\Bbb{P}^1$: 

\begin{defi} A bundle $\mathcal{V}$ on $\Bbb{P}^1$ is called generic if in any family
$\mathcal{V}_t$ of deformations of   $\mathcal{V}_0=\mathcal{V}$ over a base $T$, there is a  Zariski open neighborhood $U$ of  $t=0$, so that  $\mathcal{V}_t$ is isomorphic to $\mathcal{V}$ for $t\in U\Leftrightarrow  H^1(\pone,\mathcal{E}nd(\mv))=0\Leftrightarrow \mathcal{V}=\oplus\mathcal{O}_{\Bbb{P}^1}(a_i)$  with every $\mid a_i-a_j \mid \leq 1$. 
\end{defi}

It is easy to see that there are generic bundles of every degree and rank. Let
$\mg_{D,n}$ be the generic bundle of degree $-D$ and rank $n$.

We can therefore define
$$\langle\sigma_{I^1},\dots,\sigma_{I^{s}}\rangle_{d,D,n},$$ by replacing
in the previous definition $\mathcal{O}_{\Bbb{P}^1}^n$ by $\mg_{D,n}$.

We note the following `transformation formulas' which show that these numbers can be obtained from the usual ones.
\begin{enumerate}
\item Let $n_1,\dots,n_s$ be numbers in $\{0,\dots,n-1\}$. Let $I^1,\dots,I^s$
be subsets of $\{1,\dots,n\}$ each of cardinality $r$. Let $k_l$ be the cardinality of the set $\{u\in I^l\mid u \leq n_l\}$. 

$$\langle\sigma_{I^1},\dots,\sigma_{I^{s}}\rangle_{d,D,n}=\langle\sigma_{I^1-n_1},\dots,\sigma_{I^{s}-n_s}\rangle_{d-\sum k_l,D-\sum n_l,n}.$$ 
\item
 $$\langle\sigma_{I^1},\dots,\sigma_{I^{s}}\rangle_{d,D,n}=\langle\sigma_{I^1},\dots,\sigma_{I^{s}}\rangle_{d+r,D+n,n}.$$  
\end{enumerate}
\subsection{Grassmann Duality}\label{GDuality}
Let $W$ be a $n$-dimensional comples vector space. In this section we denote the dual of a vector space $V$ by $V^*$. We have a natural isomorphism

$$dual:Gr(r,W)\to Gr(n-r,W^{*})$$ obtained by taking $V\subset W$ to
$(W/V)^*\subset W^{*}$. 

A complete flag $F_{\bull}$ on $W$ gives rise to a complete flag $G_{\bull}$
on $W^{*}$. $G_{n-i}=$ kernel of $W^*\to F_i^*$.

\begin{lemma} Let $I=\{i_1<i_2<\dots<i_r\}\subset\{1,\dots,n\}$. Let
$$J=dual(I,n)=\{a\in \{1,\dots,n\}\mid n+1-a\not\in I\}.$$ Then

$$dual(\Omega_I(F_{\bull}))=\Omega_J(G_{\bull}).$$
\end{lemma}

\begin{proof} Let $V\in\Omega_I(F_{\bull})$ (a generic point), and $T=dual(V)$. 
$dim(T\cap G_{a})= dim(Ker (W/V)^*\to F_{n-a}^*)$
$=dim(coker(F_{n-a}\to W/V))$
$=(n-r)-dim(F_{n-a}/F_{n-a}\cap V)$
so 
$dim(T\cap G_{a})\neq dim(T\cap G_{a-1})$
if and only if 

$$n-a-dim(F_{n-a}\cap V)\neq n-a+1- dim(F_{n-a+1}\cap V)$$

that is 
$dim(F_{n-a+1}\cap V)-dim(F_{n-a}\cap V)\neq 1$
or that
$n-a+1\not\in I$.
\end{proof}

It is easy to check that $\mv_{n,D}^*$ is $\mv_{n,-D}$. Given a subbundle  $\mv $ of $\mt=\mv_{n,D}$, we get a subbundle $(\mt/\mv) ^* \subset \mt^*$ and
$degree(\mt/\mv)^*= -(deg(\mt)-deg(\mv))=D+deg(\mv)$.

\begin{lemma}
Let $I^1,\dots,I^s$ be subsets of $\{1,\dots,n\}$ each of cardinality $r$.
Let $J^l=dual(I^l,n)$ as defined in the previous lemma.

We then have 

$$\langle\sigma_{I^1},\dots,\sigma_{I^{s}}\rangle_{d,D,n}
=
\langle\sigma_{J^1},\dots,\sigma_{J^{s}}\rangle_{d-D,-D,n}$$ 

\end{lemma}

\subsection{Pre-Parabolic bundles}\label{PrePar}
\begin{defi}
Let $\mt$ be a  bundle of degree $-D$ and rank $n$ on a curve $C$ (smooth projective) endowed with $s$ complete flags on the fibers at the points $p_1,\dots,p_s$ and weights $\theta^l_a: l=1,\dots s; a=1,\dots n$ so that for any $l$

$$\theta^l_1\geq\theta^l_2\geq\dots\geq\theta^l_n\geq\theta^l_1-1.$$(No other assumptions). Call this data a pre-parabolic bundle. Also define $\theta^l_0=\theta^l_n+1$. (the pre suffix is because we have not put any absolute bounds on the $\theta^l_k$). 
\end{defi}

\begin{remark} It is usual to take partial flag varieties and strictly decreasing weights above. But we have required full flags for ease of notation. However
one can extract a partial flag variety at each parabolic point and a strictly decreasing choice of weights from each pre-parabolic bundle above. In particular
isomorphism of pre-parabolic bundles will always be by disregarding the parts of the flag varieties where the weight does not jump.
\end{remark}

Definitions of slope $\mu(\mathcal{T})$, stability and semistability can be made in the standard way. Also subbundles and quotient bundles of pre-parabolic bundles get pre-parabolic structures in the standard way. 

The following lemma is standard. 
\begin{lemma} Suppose $\mt$ is a parabolic bundle on a curve $C$ and $\mv$ a 
subbundle, $\mq$ the quotient. These acquire parabolic structures. Let $\mw\subset\mt$ be a subbundle.  Denote by $\tmw$ the saturation of the image of $\mw$ in $\mq$, and parabolic slopes by $\mu$. We then have:

$$\mu(\mw)\leq \frac{\mu(\mv)rank(\mw\cap\mv) + \mu(\tmw)rank(\tmw)}{rank(\mw)}$$
\end{lemma}

We come to the question of Harder-Narasimhan filtrations. We recall the original arguments and convince ourselves that they will go through in our context.  Let $\mathcal{T}$ be a pre-parabolic bundle on $\Bbb{P}^1$ (or for that matter any curve). If $\mathcal{T}$ is not semistable, we look for a subbundle $\mathcal{V}\subset \mathcal{T}$ satisfying
\begin{enumerate}
\item $\mathcal{V}$ is semistable.
\item For any subbundle $\mathcal{W}$ of $\mathcal{Q}=\mathcal{T}/\mathcal{V}$  we have $\mu(\mathcal{W})<\mu(\mathcal{V})$
\end{enumerate}

{\em {Existence}}: Let $\mathcal{V}$ be the subbundle of maximum slope and maximum
rank among all the subs with that slope. Existence of such a subbundle is standard.

{\em{Uniqueness}}: Let $\mv_1$ and $\mv_2$ be two such (satisfying the conditions 1,2 above. Their slopes need not be the same), look at 

$$\mv_1 \to \mt/\mv_2,$$

with image $\mw$, let $\tilde{\mw}$ be the closure of $\mw$ to a subbundle of
$\mt/\mv_2$. We have a pre-parabolic structure on $\mw$ coming from $\mv_1$,
and one on $\tilde{\mw}$ coming from $\mt/\mv_2$. Our assumptions give

\begin{enumerate}
\item $\mu(\mw)\geq  \mu(\mv_1)$.
\item $\mu(\mw)\leq \mu(\tilde{\mw})$
\item $\mu(\mv_2)> \mu(\tilde{\mw})$
\end{enumerate}

The first and second one are from the assumptions on $\mv_1$ and $\mv_2$. The one in the middle is once again okay because we are assuming the pairwise difference of $\theta^l_k$ for fixed $l$ to be $\leq 1$.

Hence $\mu(\mv_2)>\mu(\mv_1)$, and also the symmetric opposite and hence a contradiction. Hence the uniqueness.

\begin{enumerate}
\item Let $q_1,\dots q_s$ be real numbers. Let
$\tilde{\theta}^l_a =\theta^l_a +q_l$. Then this new structure ($\mt$,weights)
is (semi)-stable if and only if  the previous one is, and the Harder-Narasimhan filtrations coincide (all the slopes change by the same). This can be used to normalise:
i.e to assume

$$\sum_{a=1}^{a=n} \theta^l_a=0$$

 for each $l$.

\item The basic operation in the transformation formula is invariant under this:Let $\tilde{\mt}=$ meromorphic sections $s$ of $\mt$ so that $ts$ is a holomorphic section of $\mt$ with fiber in $F^{p_1}_1$ at $p_1$ ($t$ is a uniformising parameter there). $\tilde{\mt}$ has $s$ complete flags and weights (definition) same as those of $\mt$ except at $p_1$

$\tilde{\theta}^1_1=\theta^1_2,\dots, \tilde{\theta}^1_n=\theta^1_1-1$.

Then $\tilde{\mt}$ is semistable if and only if  $\mt$ is semistable and the Harder-Narasimhan filtrations correspond to each other under the closing up operation.
\end{enumerate} 

\subsection{The full Harder-Narasimhan Filtration}

Let $\mt$ be a generic parabolic bundle on $\Bbb{P}^1$ of degree $-D$ and rank $n$ with parabolic structure at $p_1,\dots,p_s$. Suppose it is not semistable then there exists by iterating the procedure of the previous section a filtration

$$0\subset\mv_0\subset\mv_1\subset\mv_2\subset\dots\subset\mv_k= \mt,$$

where if we denote $\mq_i=\frac{\mv_i}{\mv_{i-1}}$ we have
\begin{enumerate}
\item $\mq_i$ are semi-stable bundles.
\item $\mu(\mq_i)< \mu(\mq_{i-1})$.
\end{enumerate}

The main observation\footnote{ This is standard, in fact the stronger assertion below is an immediate consequence of the existence of the Harder-Narasimhan polygon. } is: each $\mv_i$ is unique in its ``Schubert position''.
Namely if 

$$\mv_i\in \Omega^o_{I^l}(F_{\bull}(\mt_{p_l}))$$

then
 $\langle\sigma_{I^1},\dots,\sigma_{I^{s}}\rangle_{d,D,n}=1$
where $d=-degree(\mv_i)$. 
\begin{proof}
We prove the stronger statement by induction:
if $\mv_i\neq\mw\subset\mt,rank(\mw)\geq rank(\mv_i)$, then
$\mu(\mw)<\mu(\mv_i)$ which shows the uniqueness in the Schubert class of 
$\mv_i$ immediately.

Assume the contrary, assume $\mu(\mv)=\mu(\mw)$. We will show $\mw\supset \mv_1$ and by induction we will be done. Consider the map $\mw\to\mt/\mv_1$.

$0\to\ms\to\mw\to\ml\to 0$ (with quotient structures coming from $\mw$)

We have(if $\mv_1\not\subset\mw$)
\begin{enumerate}
\item $\mu(S)<\mu(\mv_1)=\mu_1$
\item $\mu(\ml)<\mu(\mv_i/\mv_1)=\mu_2$. This is because $\ml \to \mt/\mv_1$
is a map of parabolic bundles and by induction the saturation of the image 
has less parabolic degree than $\mv_i/\mv$ (its rank is higher than that of
$\mv_i/\mv_1$).
\end{enumerate}

So $\mu(W)< \frac{rank(\ms)\mu_1 +rank(\ml)\mu_2}{rank(\ms)+rank(\ml)}$
so we will be done if we show

$$\frac{rank(\ms)\mu_1 +rank(\ml)\mu_2}{rank(\ms)+rank(\ml)}\leq\frac{rank(\mv_1)\mu_1 +rank(\mv_i/mv_1)\mu_2}{rank(\mv_1)+rank(\mv_i/\mv_1)}.$$

which simplifies to $$(rank(\mv_1)rank(\ml)-rank(\ms)rank(\mv_i/\mv_1))(\mu_1-\mu_2)\geq 0$$

But use $\mu_1>\mu_2$ and $rank(\ms)+rank(\ml)\geq rank(\mv_1)+rank(\mv_i/\mv_1)$, so 

$$(rank(\mv_1)rank(\ml)-rank(\ms)rank(\mv_i/\mv_1))\geq (rank(\mv_1)rank(\ml)$$

$$-(rank(\ms)+rank(\ml)-rank(\mv_1))rank(\ms)$$

=$$(rank(\mv_1)-rank(\ms))(rank(\mw))\geq 0$$

because $\ms\subset\mv_1$ and we are done.
\end{proof}

\subsection{Quantum Saturation}\label{QSAT}
Let us recall the geometry of the saturation problem (the reader can assume
$d=D=0$ for a first reading):

First easy case of classical saturation via a geometric argument. Start with

$$\langle\sigma_{I^1},\dots,\sigma_{I^{s}}\rangle_{d,D,n}\neq 0.$$ 

(and the case of expected dimension $=0$)
As before let $\mt=\mv_{n,D}$ with $s$ complete (and generic) flags
on the fibers of $\mt$ at $p_1,\dots,p_s$. Let $\mv$ be a point of intersection
above, let $r=rank(\mv)$. Let $\mq=\mt/\mv$. We have $H^1(\Bbb{P}^1,\mq)=0$.
We also obtain flags on the fibers of $\mv$ and $\mq$ at the parabolic points.

 The transversality at $\mv$ translates to:(and is equivalent to)

$$\overline{HOM}:=\{\phi\in Hom(\mv,\mq=\mt/\mv)\mid \phi_{p_l}(F_a(\mt_{p_l}))\subset F_{i^l_a-a}(\mq_{p_l}),$$

$$ l=1,\dots,s; a=1,\dots,r\}$$

being of the expected dimension

$$\text{expdim}=r(n-r)+dn-Dr-\sum_l(\text{codim}(\Omega^o_{I^l}(F_{\bull}(\mathcal{T}_{p_l})))=0.$$

[see ~\cite{b4} for a characterisation in terms of vector bundles for a nonvanishing GW number, $\overline{HOM}$ is just the intersection in the tangent space of the quot scheme at $\mv$ of the various inverse images of Schubert cells. The tangent space to the quot scheme is $Hom(\mv,\mq)$ from Grothendieck's theory]

$$\langle\sigma_{I^1},\dots,\sigma_{I^{s}}\rangle_{d,D,n}\neq 0$$ 
(and the case of expected dimension $=0$).
As before let $\mt=\mg_{D,n}$ with $s$ complete (and generic) flags
on the fibers of $\mt$ at $p_1,\dots,p_s$. Let $\mv$ be a point of intersection
above, let $r=rank(\mv)$. Let $\mq=\mt/\mv$. We have $H^1(\pone,\mq)=0$. Consider

$$\mt'=\mv\oplus \mq\oplus \mq\oplus\dots\oplus\mq$$

where $\mq$ appears $N$ times in the direct sum.

$\mv$ is a subbundle of $\mt'$. In the Quot scheme of subbundles of $\mt'$
, $\mv$ is a smooth point because of

$$H^1(\pone,\mathcal{H}om(\mv,\mt'/\mv))=\bigoplus_{l=1}^{l=N}H^1(\Bbb{P}^1,\mathcal{H}om(\mv,\mq))=0.$$ 

It should be noticed that $\mt'$ may not be a generic bundle. Note that
\begin{enumerate}
\item Rank$(\mt')=r+N(n-r).$
\item degree$(\mt')=-(d+N(D-d)).$
\end{enumerate}

One also gets
s partial flags on $\mq'=\mt/\mv$ by

$$F_{Na}(\mq'_{p_l})=\bigoplus_{l=1}^{l=N}F_{a}(\mq_{p_l}).$$ 

So that under each of  the $N$ natural maps $Proj_v:\mq'\to\mq; v=1,\dots,N$ one has

$$Proj_v(F_{Na}(\mq'_{p_l}))\subset F_{a}(\mq_{p_l}).$$ 

One can now induce a partial flag on the fibers of  $\mt'$ at the $s$ given points by

$$F_{N(i_l-l)+l}(\mt'_{p_l})=F_{N(i_l-l)}(\mq'_{p_l})\oplus F_{l}(\mv'_{p_l}).$$ 
for $l=1,\dots,r$.

Extend this to a complete flag on $\mt'_{p_l}$ for $l=1,\dots,s$ so that the 
Schubert position of  $\mv$ is $J^l\subset\{1,\dots,N(n-r)+r\}$, with

$$J^l=\{j^l_1<j^l_2<\dots<j^l_r\}, j^l_a=N(i^l_a-a)+a.$$

Claim: the point $\mv\subset\mt'$ defines an isolated point in the appropriate intersection, hence we find

\begin{equation}\label{eq}
\langle\sigma_{J^1},\dots,\sigma_{J^{s}}\rangle_{d,d+N(D-d),r+N(n-r)}\neq 0.
\end{equation}

The reasons:
Let $\muu$ be a universal family of  bundles of  the same degree and rank as $\mathcal{T}'$. Let $u=0$ define $\mt'$ and let $\mt'_u$  be the bundle corresponding to $u\in \muu$. Let $Quot$ be the scheme over $\muu$ of  Quotients of  $\mathcal{T}'_u$ of  the same rank and degree as $\mq'$. It is standard that $Quot$ is projective over $\muu$. We have 
\begin{enumerate}
\item The subbundle $\mv \subset\mt'$ defines a smooth point $q_0$ in  the fiber of  $Quot$ over $u=0$. The map $Quot\to \muu$ is hence smooth at this point. The map is flat because the vanishing of   $H^1(\pone,\mathcal{H}om(\mv,\mt'/\mv))$  removes all obstructions to motion (the obstructions all live in $H^1(\pone,\mathcal{H}om(\mv,\mt'/\mv))$ (see \cite{kollar} page 31, 2.10.4 for an analogous assertion).
\item $\overline{HOM}$ corresponding to ~\ref{eq} is a $N$-fold direct sum of 
the $\overline{HOM}$ corresponding to 

$$\langle\sigma_{I^1},\dots,\sigma_{I^{s}}\rangle_{d,D,n}\neq 0.$$ 

when evaluated at $\mv\subset \mt'$ with the above flags on $\mt'$ and
at $\mv \subset \mt$ respectively. That is if,
$$\overline{HOM}_1:=\{\phi\in Hom(\mv,\mq=\mt/\mv)\mid \phi_{p_l}(F_a(\mt_{p_l}))\subset F_{i^l_a-a}(\mq_{p_l}),$$

$$ l=1,\dots,s; a=1,\dots,r\}$$

then
$$\overline{HOM}_2:=\{\phi\in Hom(\mv,\mq'=\mt'/\mv)\mid \phi_{p_l}(F_a(\mt_{p_l}))\subset F_{N(i^l_a-a)}(\mq'_{p_l}),$$

$$ l=1,\dots,s; a=1,\dots,r\}$$

satisfies

$$\overline{HOM}_2=\bigoplus_{a=1}^{N}\overline{HOM}_1.$$

But our assumption of $\mv$ being a transverse point of intersection in the intersection defining 
$$\langle\sigma_{I^1},\dots,\sigma_{I^{s}}\rangle_{d,D,n}\neq 0$$

tells us that $\overline{HOM}_1=\{0\}$, so this gives $\overline{HOM}_2=\{0\}$ or that $\mv$ is a transverse point of intersection (at a smooth point) of the Schubert cells defining $\langle\sigma_{J^1},\dots,\sigma_{J^{s}}\rangle_{d,d+N(D-d),r+N(n-r)}$ when we substitute \\$\mg_{d+D(D-d),r+N(n-r)}$ by $\mt'$ (which has the same rank and degree) and the flags constructed above on $\mt'$. From local intersection theory we then get $$\langle\sigma_{J^1},\dots,\sigma_{J^{s}}\rangle_{d,d+N(D-d),r+N(n-r)}\neq 0.$$

\end{enumerate}

Introduce the notation: for $I$ a subset $\{i_1<\dots\i_r\}$ of $\{1,\dots,n\}$, define $NI=\{k_1<\dots<k_r\}\subset \{1,\dots,r+N(n-r)\}$ by $k_l= l+N(i_l-l),l=1,\dots r$.

Hence we can form the function:
\begin{equation}\label{deff}
f(N)=\langle\sigma_{NI^1},\dots,\sigma_{NI^{s}}\rangle_{d,d+N(D-d),r+N(n-r)}.
\end{equation}
We have shown the easy half of the Saturation problem:

$f(1)>0\Rightarrow f(N)>0$, for all $N$.

The proof of the other half is in ~\cite{b4}. 
\subsection{Eigenvalue form of Quantum Horn conjecture}

Consider a situation where
 
$$\langle\sigma_{I^1},\dots,\sigma_{I^{s}}\rangle_{d,D,n}\neq 0.$$ 

Let generic flags be chosen on the fibers of $\mt=\mg_{D,n}$ at the points $p_l,l=1,\dots,s$. And let $\mathcal{V}$  be a point on the intersection `defining' 
$$\langle\sigma_{I^1},\dots,\sigma_{I^{s}}\rangle_{d,D,n}\neq 0.$$

\begin{remark}(Genericity Property:) Note that we can assume that $\mv$ is generic as a bundle (together with Kleiman Bertini theorem).  
$\mathcal{V}$ inherits complete flags on its fibers, and can be assumed to be isomorphic
to $\mg_{d,r}$ and the flags generic enough for intersection theory
(the reason for this is same as in \cite{b3},  the idea is the following: any small
perturbation of the flags on $\mg_{D,n}$ is still generic, perturb it is such a way that $\mv$ stays in the intersection and the induced flags on it
become generic). 
\end{remark}

 Define numbers
$\alpha^l_k=\frac{n-r+k-i^l_k}{n-r}$ and these are numbers between $0$ and $1$
(both possible). This gives the structure of a pre-parabolic bundle on $\mathcal{V}$. It is the easy half of the (Quantum Analogue) of the Horn Conjecture ~\cite{b4} , that this structure is semistable of slope $1+\frac{D-d}{n-r}$. The transformation formulas give various `translates' of this pre-parabolic bundle. One can always assume that the upper limit
$1$ is never attained by the $\alpha$'s. Also one can assume $D=0$, but without the $\alpha$ assumption.

\begin{defi} Let $\mathcal{T}$ be a bundle of degree $-D$ and rank $n$ on $\Bbb{P}^1$, endowed with flags $F_{\bull}(\mathcal{T}_{p_l}), l=1,\dots,s$. Let $\mathcal{V}$
be a subbundle of $\mathcal{T}$ of rank $r$ degree $-d$ whose fiber at the point $p_l$ is in the Schubert cell $\Omega^o_{I^l}(F_{\bull}(\mathcal{T}_{p_l}))$. Define

$$dim(\mathcal{V},\mathcal{T}, \{F_{\bull}(\mathcal{T}_{p_l}), l=1,\dots,s\})$$

$$=r(n-r)+dn-Dr-\sum_l(\text{codim}(\Omega^o_{I^l}(F_{\bull}(\mathcal{T}_{p_l}))).$$

That is the expected dimension of subbundles of $\mt$ in the same `Schubert Position' as $\mathcal{V}$. Call the data $(d,r,I^1,\dots,I^s)$ the Schubert position of $\mathcal{V}$. When the flags on $\mathcal{T}$ are clear we will just write
$dim(\mathcal{V},\mathcal{T})$ for the dimension above.
\end{defi}

\begin{remark}\label{stmnt} The semistability inequalities for a subbundle $\ms\subset\mv$
are the same as
$$dim(\ms,\mv,\{F_{\bull}(\mathcal{V}_{p_l}), l=1,\dots,s\}) \leq dim(\ms,\mt,\{F_{\bull}(\mathcal{T}_{p_l}), l=1,\dots,s\})$$

where $F_{\bull}(\mathcal{V}_{p_l})$ are the induced flags on the fibers of $\mv$. (see ~\cite{b4} for details)
\end{remark}
\begin{defi} The pre-parabolic bundle $\mv$ defined above is called the Witten bundle corresponding to the intersection defining 
$$\langle\sigma_{I^1},\dots,\sigma_{I^{s}}\rangle_{d,D,n}\neq 0.$$ 

To be precise, this pre-parabolic bundle depends on the choice of the flags 
used to evaluate the intersection (and the point chosen). This bundle appears
in a theorem of Witten ~\cite{witten} on Quantum cohomology as will be explained later in section ~\ref{HWitten}. 
\end{defi}
\begin{remark}
We need to generalize this construction a bit to handle the following situation: Let $\mt$ be a parabolic whose underlying bundle is generic, and $p$ is a parabolic point. It is important to disregard the parts of the flag $F(\mt_p)_{\bull}$ where the weights do not jump. 

So assume we are given a generic bundle of degree $-D$ and rank $n$ and for each $p\in\{p_1,\dots,p_s\}$ we are given a flag:

$$F(\mt_p)_{u^p_1}\subset\dots F(\mt_p)_{u^p_{k(p)}}=F(\mt_p)$$

where the elements on the subscripts are the dimensions. Assume that these partial flag varieties are in general position. Also assume that we have particular extensions of these flags into complete flag varieties. Suppose also that $d$ and $r<n$ are given. Now given any string of numbers $v^p_1\leq v^p_2\leq \dots\leq v^p_{k(p)}=r$ we can define a Schubert cell\footnote{The closure of this cell is a standard Schubert cell when the flags of $\mt_p$ are extended to a full flag.} in $Gr(r,\mt_p)$ by 

$$\{V\subset \mt_p\mid dim(V\cap F(\mt_p)_{u^p_{l}}=v^p_l, l=1,\dots k(p)\}$$

Denote the closure by $\Omega_{I^l}(\mt_p)$, where the data $I^l$ can be determined from the earlier data. Suppose that 

$$\langle\sigma_{I^1},\dots,\sigma_{I^{s}}\rangle_{d,D,n}\neq 0.$$ 

That is there is a subbundle $\mv$ of $\mt$ of degree $-d$, rank $r$ whose fiber at each $p\in\{p_1,\dots,p_s\}$ satisfies

$$\{V\subset \mt_p\mid dim(V\cap F(\mt_p)_{u^p_{l}}=v^p_l, l=1,\dots k(p)\}$$

and there are only finitely many such $\mv$ (we have used Kleiman Bertini theorem).

From before $\mv$ as a Witten bundle acquires weights. We note that the weights
at $p$ satisfy: weight of $F(\mv_p)_{v^p_l}=\frac{n-r+v^p_l-u^p_l}{n-r}$. If $v\not\in\{v^p_1,\dots,v^p_k(p)\}$,and $v_a$ is the smallest element in $\{v^p_1,\dots,v^p_k(p)\}$ satisfying $w<v_a$, then weight of $F(\mv_p)_w=\frac{n-r+v^p_a-u^p_a}{n-r}$. So the weights of the Witten bundle can increase only at $\mv  
\cap F(\mt)_{u^p_l}$. So the extension of the partial flag on $\mt_p$ to a full one is irrelevant.
\end{remark}
It is natural to consider the moduli space of pre-parabolic bundles of the type of $\mv$. But the Geometric invariant theory method ~\cite{Pauly} requires
$\alpha^l_1<\alpha^l_r+1$ (that is successive differences of $1$ are not allowed). The GIT quotient is then the natural equivalence (Jordan-Holder equivalence). Without assumptions on $\alpha$ we do not have an abelian category of semistable bundles.

Using the transformation formulas we can assume that the difference between the $\alpha's$ at any point are strictly less than $1$. It is useful to make the following definition.

\begin{definition}\label{NORMAL} We will call a cycle $I$ normalised if this is the case. This is equivalent to either $i_1>1$ or $i_r<n$.
\end{definition}

It should be noted that for any two choices of the correction the Moduli Spaces are in canonical correspondence\footnote{With the $\Theta$ bundles corresponding under this isomorphism.} (with the equivalence, and depends on the choice of the correction). We will call this the Witten Moduli Space of $\mv$. The Witten bundle gives a bundle in this moduli space (but there is more identification if the successive difference of some $\alpha$'s at a point is $1$).

\begin{example} $Gr(2,4)$,$s=2m$, $I^l=\{1,4\}, l=1,\dots s$, $d=m-1$. The associated Witten Moduli Space is a point.
\end{example}

\begin{defi} If the pairwise difference of the $\alpha$'s is strictly less than $1$  we will call the cycle normalised.
\end{defi}

\begin{remark}
We will find it convenient to use expressions of the form 'the Witten bundle is stable' - to mean the moduli space has a stable point (a nonempty open set). Or, ``the Witten bundle is rigid'' to mean the dimension of the Witten Moduli space $=0$. \end{remark}

We now state the Eigenvalue form  of the quantum Horn Problem.

\begin{defi}Let 

$$\Delta_r =\{(a_1,a_2,\dots,a_r)\mid 
a_1\geq a_2\dots\geq a_r\geq a_1-1, \sum a_i=0\}.$$
\end{defi}

\begin{defi}
For ${\bf {a}}=(a_1,\dots,a_r)\in \Delta_r$, define

$$S({\bf{a}})=(a_2+\frac{1}{r},\dots,a_r+\frac{1}{r},a_1-(1-\frac{1}{r})).$$
\end{defi}
Also define

\begin{defi}$I=\{ i_1< \dots< i_r \}$ subset of  $\{1,\dots,n\}$ of cardinality $r$
Define $\delta(I)=(\delta_1,\dots,\delta_r)\in \Delta_r$, by

$\delta_j =\frac{j-i_j}{n-r} -c$
where $c = \sum_{j=1}^{r}\frac{j-i_j}{(r(n-r))}$.
\end{defi}
We have therefore proved the easy half of
\begin{theorem}\label{Evalues}(The Quantum Horn Conjecture)Let $I^1,\dots,I^s$ be subsets of
$\{1,\dots,n\}$ each of cardinality $r$. Use the notation

$$I^j=\{i^j_1<i^j_2<\dots<i^j_r\}.$$ Assume that

$$\sum_{j=1}^{j=s}\sum_{a=1}^{a=r} (n-r +a -i^j_a)= r(n-r)+nd.$$
where d is a non negative integer
then 
$$<\sigma(I^1),\dots,\sigma(I^s)>_d \neq 0$$

if and only if:

                Given any $1 \leq p < r$ and any choice of subsets $K^{1},K^{2},\dots,K^{s}$ of cardinality $p$ of $\{1,\dots,r\}$ and $q$ a nonegative integer and if 
$$<\sigma(K^1),\dots,\sigma(K^s)>_q \neq 0$$

then the  inequality

$$\sum_{a\in K^{1}} (S^d(\delta(I^1))_a + \sum_{j=2}^s \sum_{a\in K^{j}}(\delta(I^j)_a) \leq q$$ 

is valid.

\end{theorem} 

\begin{remark} The operation $S^d$ is important. 
Let us look at $Gr(2,4)$. Here $r=2,n=4,d=1$, 
$I^1,I^2,I^3$ are $\{1,4\},\{2,3\}$ and $\{1,2\}$, the corresponding $\delta's$ are
$(1/2,-1/2),(0,0)$ and $(0,0)$, and the shift operation on the first converts it
into $(0,0),(0,0),(0,0)$ which certainly satisfies the inequalities. But without the shift, the inequality corresponding to
$K^1,K^2,K^3=\{1\},\{2\},\{2\}, p=1,q=0,$ would have led us to
$\frac{1}{2}\leq 0$ which is false.

\end{remark}

\begin{remark}\label{remarkable} Choose $d_1,\dots,d_s$ non negative with sum $=d$ (and fix them once and for all), then the inequalities could be written as

$$ \sum_{j=1}^s \sum_{a\in K^{j}}S^{d_j}(\delta(I^j))_a \leq q.$$
 \end{remark}
\begin{remark} In the formulation of remark ~\ref{remarkable} the conditions are the same as there to exist matrices $A_l, l=1,\dots,s$ in $SU(r)$ with conjugacy class of $A_l$ being $S^d_l(\delta(I^l)), l=1,\dots,s$ and $A_1A_2\dots  A_s=1$ (see \cite{AW}, \cite{b1}).
\end{remark}

\subsection{The theorem of Witten on Quantum Cohomology}\label{HWitten}
The principal aim of this section is to justify the following claim. Consider the situation when we are considering a GW number (in the case expected dimension =0)
$$\langle\sigma_{I^1},\dots,\sigma_{I^{s}}\rangle_{d,D,n}$$
where the cycles are normalised. Form the function $f(N)$  as in definition ~\ref{deff}. We claim the equivalence of the two statements
\begin{enumerate}
\item $f(N)=1,\forall N>0.$
\item The Witten moduli space $\mathcal{M}$ of $\langle\sigma_{I^1},\dots,\sigma_{I^{s}}\rangle_{d,D,n}$ is a point
\end{enumerate}

Obviously if $f(N)=h^0(\mathcal{M},\Theta^N)$ for $\Theta$ an ample line bundle on $\mathcal{M}$, then this claim is immediate. For this we need a mechanism
of writing GW numbers as global sections of ample line bundles on suitable moduli spaces. This is available by the work of Witten ~\cite{witten}.

In ~\cite{witten}, Gromov-Witten invariants of Grassmannians are related to
dimensions of global section of line bundles over suitable moduli spaces of parabolic bundles.
That is, an equation of the form 
$$\langle\sigma_{I^1},\dots,\sigma_{I^{s}}\rangle_{d,D=0,n}=h^0(\mathcal{M},\theta)$$
is proved. We need to describe the modulispace $\mathcal{M}$ and the line bundle $\Theta$ over it.  The moduli space is the moduli space of representations $\mathcal{M}$ of $\pi_1(\Bbb{P}^1-\{p_1,\dots,p_s\})$ into $SU(r)$ with monodromy at $p_i$ given by
$S^{d_j}(\delta(I^j))$ as above (where $\sum d_j=d$ as in lemma ~\ref{remarkable}). This has a natural level structure ($n-r$) and hence as in ~\cite{Pauly} a line bundle $\Theta$ on it. This gives the $\mathcal{M}$ and $\Theta$ in Witten's theorem.  

Using the transformation formulas, one sees now that Generalised GW numbers of the form $\langle\sigma_{I^1},\dots,\sigma_{I^{s}}\rangle_{d=0,D,n}$ are $h^0(\mathcal{M},\Theta)$ where $\mathcal{M}$ is the moduli space of representations $\mathcal{M}$ of $\pi_1(\Bbb{P}^1-\{p_1,\dots,p_s\})$ into $SU(r)$ with monodromy at $p_i$ given by $\delta(I^j)$(no shifting) . This has a natural level structure ($n-r$) and hence as before a line bundle $\Theta$ on it. 

In general $\langle\sigma_{I^1},\dots,\sigma_{I^{s}}\rangle_{d,D,n}$ are $h^0(\mathcal{M},\Theta)$ where $\mathcal{M}$ is the moduli space of representations $\mathcal{M}$ of $\pi_1(\Bbb{P}^1-\{p_1,\dots,p_s\})$ into $SU(r)$ with monodromy at $p_i$ given by $S^{d_j}(\delta(I^j))$ as above (where $\sum d_j=d$ as in lemma ~\ref{remarkable}). This has a natural level structure ($n-r$) and hence as before a line bundle $\Theta$ on it. This parabolic modulispace is isomorphic to the Witten modulispace of $\langle\sigma_{I^1},\dots,\sigma_{I^{s}}\rangle_{d,D,n}$ in the case $I^j$ are normalised.

Let $\mathcal{M}_1$ be the Witten moduli space of pre-parabolic bundles corresponding to $\langle\sigma_{I^1},\dots,\sigma_{I^{s}}\rangle_{d,D,n}\neq 0.$ The
equivalence imposed is isomorphism of the Jordan-Holder factors (which are stable)\footnote{if both limits are allowed for the $\alpha$'s, we no longer get an Abelian category of semistable bundles.}. By Grassmann duality we get a similar moduli space of Quotients $\mathcal{M}_2$. Each of these have natural ample line bundles $\Theta_1$ and $\Theta_2$, and the theorem of Witten gives \footnote{all the non standard GW numbers can be reduced to the usual GW numbers.}

$$\langle\sigma_{I^1},\dots,\sigma_{I^{s}}\rangle_{d,D,n}=h^0(\mathcal{M}_1,\Theta_1)=h^0(\mathcal{M}_2,\Theta_2).$$

Also the function
$$f(N)=\langle\sigma_{NI^1},\dots,\sigma_{NI^{s}}\rangle_{d,d+N(D-d),r+N(n-r)}$$
from the earlier section\footnote{One should make explicit what $\Theta_1$ is to do so. Naively the local factors of $\Theta_1$ come from the map of $\mathcal{M}_1$ to the partial flag varieties and there is global twisting by determinant of cohomology and by determinant of a general fiber. This expression can be found for example in ~\cite{Pauly}. In $f(N)$ we are scaling the local factors by $N$ and we need to make sure that the global corrections match too.}

$$f(N)=h^0(\mathcal{M}_1,\Theta_1^N).$$

The quantum generalization of Fulton's conjecture is the following:
\begin{question}
Does
$f(1)=1\Rightarrow f(N)=1$, for any $N$?
\end{question} 

[The other implication is true by the Quantum Saturation conjecture, and Witten's theorem (which gives $f(1)\leq f(N)$).]

The link to the Witten moduli space is via the observation (because of the ampleness of $\Theta_1$): 

$$\forall N>0, f(N)=1 \Leftrightarrow \mathcal{M}_1=\text{ a point. }$$

The classical part of IRQ (Fulton's conjecture)  was proved by Knutson, Tao and Woodward ~\cite{KTW}. A Geometric approach  would be to prove the base point freeness of $\mathcal{\theta}_1$ on  $\mathcal{M}_1$. This is false even in the classical case ($d=D=0$). For this would imply that $f(1)
\geq dim(\mathcal{M}_1)+1$ (``the theta functions'' will give a finite map to projective space\footnote{ If $X$ is a projective variety and $L$ an ample line bundle, $\theta_i$ global sections of $L$ which do not have any common zeroes, then the corresponding map of $X$ to projective space is finite.} ) which is false (an example will be given in example ~\ref{egg}). The existence of such examples also puts an end to the hope that ``theta functions'' of level $1$ give a generic (local analytic) embedding. 

It could very well be case that the Quantum analogue is false (we have not done any computations beyond the classical examples). The numbers involved in a counterexample could be large. The counterexample to the base point `nonfreeness' example takes place in a comparatively large Grassmannian $Gr(8,12)$.

\begin{remark}\label{decision}
We can also produce  a number in the  Quantum case $M$ so that $f(M)=1$ implies $\mathcal{M}_1=$ a point. In fact $M=r+1$ is probably sufficient. In the classical case $M=1$ is sufficient by the result of ~\cite{KTW}. The proof of the bound($M=r+1$) will be written up elsewhere. From a computational point of view
there is another way of determining whether the moduli space $\mathcal{M}$ is
a point. If it has a stable point, then there is an exact formula for the dimension of $\mathcal{M}$ (in section ~\ref{irred}). If it is not stable then we have to 
check the rigidity of each Jordan-H\"{o}lder factor. So we need to decide whether a generic point in $\mathcal{M}_1$ is semistable but not stable. This means that equality holds for some subbundle in the Quantum Horn conjecture (theorem).
One can look at all subbundles, and check if equality holds and if it does
repeat the same procedure for the subbundle and for the quotient. The subbundle  can be assumed to be  generic (with generic flags).
\end{remark}
\begin{remark} We make a remark on the Grassmann dual Witten bundle. Let
$$\langle\sigma_{I^1},\dots,\sigma_{I^{s}}\rangle_{d,D,n}\neq 0$$
and let $\mv\subset\mg_{D,n}$ be a point of intersection and let $\mq=\mg_{D,n}/\mv$. $\mq^*$ gets a semistable parabolic structure (by Grassmann duality), and hence $\mq$ gets one too. Let $\mw$ be a subbundle of $\mq$ which lifts to a subbundle $\tilde{\mw}\subset \mg_{D,n}$. The semistability inequality for $\mw$ is same as

$$dim(\mw,\mq)\leq dim(\tilde{\mw},\mg_{D,n}).$$

(all with the naturally induced flags). This statement is obtained by dualizing remark ~\ref{stmnt}.
\end{remark}

\subsection{Eigenvalue problems}\label{Evalue1}

In ~\cite{AW} and ~\cite{b1} the following problem was considered. Let
$\bar{A}_1,\dots,\bar{A}_s$, be conjugacy classes in $SU(n)$, find conditions on these so that there are elements $\bar{A}_i\in SU(n)$ with conjugacy class of
$A_i$ in $\bar{A}_i$ and

$$A_1A_2\dots A_s=1$$

This problem is  in a natural way related to unitary representations
of the fundamental group of $\Bbb{P}^1-\{p_1,\dots,p_s\}$ where $\{p_1,\dots,p_s\}$ is a set of distinct points on $\Bbb{P}^1$. The answer to this question involves Quantum Cohomology of Grassmannians.  
Recall the set up:
\begin{defi} Let $\Delta_{n-1} =\{\vec{a}=(a_1,\dots,a_n)\mid a_1\geq a_2\geq\dots\geq a_n\geq a_1-1), \sum a_l =0\}$. In a natural way, these parametrise the conjugacy classes in $SU(n)$, the conjugacy class corresponding to $\vec{a}$
is the diagonal matrix with entries $e^{2\pi i a_l}$. Let $\Delta^o_{n-1}$ be the interior of this simplex: $\Delta^o_{n-1} =\{\vec{a}=(a_1,\dots,a_n)\mid a_1 > a_2>\dots> a_n> a_1-1), \sum a_l =0\}$. 
\end{defi}

Also recall the notation:
\begin{defi}Let $\bar{A}=\vec{a}\in \Delta_{n-1}$, and $I=\{i_1<i_2<\dots<i_r\}
\subset \{1,\dots,n\}$. Define

$$\lambda_I(\bar{A})=\sum_{t\in I}a_t.$$
\end{defi}

Recall the principal result in ~\cite{b1}
\begin{theorem}Let $\bar{A}_l,l=1,\dots,s$ be conjugacy classes in $SU(n)$.
. Then there exist $A_l\in SU(n)$with conjugacy class of $A_i$ in $\bar{A}_i$ and $A_1A_2\dots A_s=1$ if and only if  
given any integers $1\leq r<n, d\geq 0$,subsets $I^1,\dots I^s$ of $\{1,\dots,n\}$ each of cardinality $r$, such that
$$\langle\sigma_{I^1},\dots,\sigma_{I^{s}}\rangle_{d}=1,$$

the inequality

$$\sum_{l=1}^{l=s}\lambda_{I^l}(\bar{A}_l) -d\leq 0$$

holds. 
\end{theorem}
Let us name the inequality corresponding to the data  $(d,r,I^1,\dots,I^s)$
as\linebreak $Ineq(d,r,I^1,\dots,I^s)$. So according to the theorem above the set of inequalities $Ineq(d,r,I^1,\dots,I^s)$ satisfying
$$\langle\sigma_{I^1},\dots,\sigma_{I^{s}}\rangle_{d}=1$$

form a necessary and sufficient set of inequalities defining the desired
set inside $\Delta_{n-1}^s$.

The reduction to the set of inequalities corresponding to intersection $=1$
was made using the existence of Harder-Narasimhan filtration. The question
pops up if there are any more such conditions. The theorem has been extended
to other groups in many contexts (Quantum, Classical), where it is either
expected or known that some of the inequalities (corresponding to intersection
 $=1$) are redundant. 

The same methods give the following theorem where as before we assume that the pairwise difference of any two $\theta$ at any point is $\leq 1$. The advantage
is that there is considerably more freedom to try to prove the irredundancy.
It should also be noted that the following has exactly the same content as the previous theorem. We can use the transformation formulas to reduce to $D=0$ and norm the $\theta$'s (see lemma ~\ref{NORM}) to sum to zero at any point.

\begin{theorem}\label{SS}
Let the weights $\theta^l_j$,$D$ and $n$ be given, then $\mt=\mathcal{G}_{n,D}$ with $s$ generic flags and
weights $\theta^l_a$ is semistable if and only if whenever given $r$, $s$ subsets  $I^1,\dots,I^s$ of $\{1,\dots,n\}$ each of cardinality $r$ satisfying

$$\langle\sigma_{I^1},\dots,\sigma_{I^{s}}\rangle_{d,D,n}=1,$$
one has

$$\frac{1}{r}(-d+\sum_l\sum_{a\in I^l} \theta^l_a )\leq$$

$$\frac{1}{n}(-D+\sum_l\sum_{a=1}^{a=n} \theta^l_a )$$
\end{theorem}

Label the inequalities by $Ineq(I^1,I^2,\dots,I^s,d,D,n).$

\begin{defi}Call an $Ineq(r,d,D,I_1,\dots,I_s,n)$ polyrigid if  in the function
$$f(N)=\langle\sigma_{NI^1},\dots,\sigma_{NI^{s}}\rangle_{d,d(1-N) +ND,r+N(n-r)},$$ 
$f(N)=1$ for all $N$ (or that the corresponding Witten Moduli Space is a point\footnote{ The main point is that $f(N)=H^0(\mathcal{M},\theta^N)$ where $\mathcal{M}$ is a moduli space of parabolic bundles (which identifies two parabolic bundles with the same Jordan-Holder factors ) and $\theta$ is ample on $\mathcal{M}$. So $f(N)=1$ for all $N$ implies that the moduli  space is a point. Hence all the graded factors must be rigid. 
}). 
\end{defi}
\subsection{What does irredundancy mean?}\label{RED1}
There are a few `different' converses to ~\ref{SS}. They are all of the form: Find a subset of the above list of inequalities (call the subset $\mathcal{I}$) so that they describe the same polyhedron and, 

\begin{enumerate}
\item Given $Ineq(d,D,I^1,\dots,I^s,n)$ in $\mathcal{I}$, there are parabolic weights so that the subbundle  corresponding to this inequality for generic choice of flags) is the Harder-Narasimhan element.\footnote {One could ask a similar question on all of the Harder-Narasimhan filtration}.

\item For each inequality $Ineq(d,D,I^1,\dots,I^s,n)$ in the list $\mathcal{I}$, there are choice of weights so that
for the corresponding parabolic bundle, the subbundle corresponding to 
$\langle\sigma_{I^1},\dots,\sigma_{I^{s}}\rangle_{d,D,n}$
is the {\em{sole}} contradiction to semistability (among all subbundles).
\end{enumerate}

\begin{remark}
Consider a non semistable parabolic bundle. It has a unique Harder-Narasimhan
subbundle  (the first part of the filtration). This subbundle  contradicts semistability, but is not necessarily the only subbundle  contradicting semistability. Example $\mv_1\oplus\mv_2\oplus\mv_3$ and assume $\mv_i$ semistable and the slope $\mu(\mv_1)>\mu(\mv_2)
>\mu(\mv_3)$. Then clearly $\mv_1$ is the Harder Narasimhan element, but
$\mv_1\oplus\mv_2$ might contradict semistability too! (if the slope of $\mv_3$
is very small).

\end{remark}

 The strongest converse is (no.2). But the most natural converse is the first one. Let us call an inequality $Ineq(d,D,I^1,\dots,I^s,n)$ HN-irredundant if the first possibility occurs, and to be strongly irredundant if the second holds. Clearly strongly irredundant $\Rightarrow$ HN-irredundant. 

\begin{defi} The large weight space $\Omega_n$ is the set:

$$\{(t_1,\dots,t_n)\mid 1\geq t_1\geq t_2\geq\dots\geq t_n\geq 0\}.$$

The interior of $\Omega_n$ is the set $\Omega^o_n$: 

$$\{(t_1,\dots,t_n)\mid 1> t_1> t_2>\dots> t_n> 0\}.$$
\end{defi}

\begin{lemma}\label{NORM} There is a natural surjection $$\Omega_n\to \Delta_n$$ which
restricts to a surjection of the interiors: $\Omega^o_n\to \Delta^o_n$.
\end{lemma}
\begin{proof} The map is $(t_i)$ to $t_i-c$ where $c=\frac{\sum t_i}{n}$.
\end{proof}

\begin{defi} Let $\Lambda_{D,n}\subset (\Omega^o_n)^s$ be the set of weights which give generically a semistable parabolic structure on $\mg_{D,n}$ with parabolic structure at $s$ points. This set has a nonempty\footnote{ There do exist irreducible unitary representations of $\pi_1(\Bbb{P}^1-\{p_1,\dots,p_s\}$ for $s\geq 3$ - for example take $A,B$ unitary matrices so that the group they generate is Zariski dense in $U(n)$ and take the representation which sends the fundamental loops (for $s=3$) to $A,B, (AB)^{-1}$}  interior  and is a polyhedron.
\end{defi}

\begin{lemma} \label{lemmacondition} An inequality $Ineq(d,D,I^1,\dots,I^s,n)$ is strongly irredundant if and only if 
\begin{enumerate}
\item There are choices of weights for which $Ineq(d,D,I^1,\dots,I^s,n)$ is not valid.
\item \label{theother} There are weights (for each point) in $\Omega_n$ which give a semistable structure on $\mathcal{G}_{D,n}$ for which (having chosen generic flags on $\mg_{D,n}$), the subbundle  corresponding to $Ineq(d,D,I^1,\dots,I^s,n)$ has the same parabolic slope as the whole space and the weights for the subbundle are separated away from the quotient. To be more specific, for a parabolic point (say $p_1$),
let $1\geq\theta_1\geq\dots\geq\theta_n\geq 0$ be the weights at $p_l$. then 
we have for $i\in I^l,j\in\{1,\dots,n\}-I^l$, 

$$0<\mid\theta_i-\theta_j\mid<1.$$
\end{enumerate}
\end{lemma}
\begin{proof} 
First for generic flags on $\mt=\mg_{D,n}$, let $\mv$ be the point of intersection corresponding to $Ineq(d,D,I^1,\dots,I^s,n)$. 

One direction is trivial. If an inequality is strongly irredundant, find a choice of weights for which $Ineq(d,D,I^1,\dots,I^s)$ is the only contradictor of semistability. Also find another choice of weights that gives a stable structure
on $\mt$ (in the interior of the large weight space). Join these two choices of weights and find the point on the line joining these weights where the structure becomes semistable (and slope of $\mv$ equals the slope of $\mt$).

 Let us look at the ``issues'' with the other direction: We have weights that make $\mg_{D,n}$ semistable, $\mv$ is a contradictor to stability  and the weights of $\mv$ are separated from those of $\mq$. We want weights for which the subbundle $\mv$ is the sole contradiction to semistability. We do this in two steps.

First step: Perturb the given weights on $\mv$ and on $\mq$ so that both become stable (there is space to do this and induce these weights to $\mg_{D,n}$). This perturbation can change the slopes of these (we are working in the large weight space) and these slopes may not be equal. We can cure this problem by adding a constant to all the weight of (say the ones corresponding to the subbundle ) at one point. This might take us out of the bounds :$[0,1]$ but the pairwise difference of weights at any point is still $<1$. So we can correct this last problem by a suitable shift of all the weights at a point(by the same amount). Hence we get a choice of weights for which $\mv$ is the only contradiction to stability (and $\mq$ if $\mv$ splits off). Therefore for any small perturbation of these weights, the only possible contradictor to stability (or semistability) is $\mv$ (and possibly $\mq$).

Second step:  Now join this point to the choice of weights in the first condition. These are choices of weights for which $\mv$ contradicts semistability and ``head out just a little''. Even if $\mv$ splits off, $\mq$ has now slope less than that of $\mt$.
\end{proof}

\section{Redundancy in Eigenvalue Problems.}\label{redundancy}

\begin{theorem}\label{SS2} In theorem(~\ref{SS}) the set any  inequality 
$Ineq(r,d,D,I_1,\dots,I_s,n)$ which is not polyrigid or whose Grassmann dual is not polyrigid  does not get to be  the Harder-Narasimhan element for any choice of weights.

Hence the set of inequalities $Ineq(r,d,D,I_1,\dots,I_s,n)$ which are polyrigid
and whose Grassmann duals are polyrigid form a necessary and sufficient set of inequalities.

\end{theorem}
\begin{proof}
 
Suppose $Ineq(r,d,D,I^1,\dots,I^s,n)$ is a non polyrigid inequality and yet there
is a choice of weights $\theta^l_i$ so that the subbundle $\mv$ (choose flags etc) corresponding to $\langle\sigma_{I^1},\dots,\sigma_{I^{s}}\rangle_{d,D,n}=1,$
is the Harder-Narasimhan element. So we have:

$$0\to \mathcal{V} \to \mt\to Q\to 0$$

and every subbundle $\mathcal{W}$ of $Q$ satisfies $\mu(\mathcal{W})<\mu(\mathcal{V})$ and $\mathcal{V}$ is semistable.

Suppose in the function $f(N)>1$ for some $N>1$ \footnote{ Recall the definition of $f(N)$  from ~\ref{deff}.}. 

Consider $\mt'=\mq'\oplus \mv$ where  $\mq'=Q\oplus\dots\oplus Q$(repeated $N$ times). One can consider $\mt'$ as a pre-parabolic bundle in the standard way (one may need to complete the flag but the weights assigned are the smallest possible to such an addition). Also note that the Schubert position of $\mv\subset\mt'$ is $J^l\subset\{1,\dots,N(n-r)+r\}$, with

$$NI^l=J^l=\{j^l_1<j^l_2<\dots<j^l_r\}, j^l_a=N(i^l_a-a)+a.$$

And as before (see section ~\ref{QSAT}) the point $\mv\subset\mt'$ defines an isolated transverse point in the appropriate intersection, corresponding to
$$\langle\sigma_{NI^1},\dots,\sigma_{NI^{s}}\rangle_{d,d+N(D-d),r+N(n-r)}\neq 0.$$

and one has the exact sequence

$0\to \mathcal{V}\to \mt'\to \mq'\to 0$. Clearly $\mathcal{V}$ is the Harder-Narasimhan element in $\mt'$, but two points

\begin{enumerate}
\item $f(N)>1$(this is our assumption).
\item $\mt'$ is not generic, but a limit of generic bundles (with generic flags). On each of them we will find $\mathcal{V}_t,\mathcal{W}_t$ so that $\mathcal{V}_t$ tends to $\mathcal{V}$
and $\mathcal{W}_t$ stays away from $\mathcal{V}$ (we are assuming $f(N)>1$) in the appropriate Quot scheme. The limit of $\mathcal{W}_t$ is a coherent subsheaf of $\mt$ and its closure is a subbundle  with at least as much parabolic degree  (and of the same rank) as $\mathcal{V}$. Hence we find a contradiction.
\end{enumerate}
\begin{remark}
If $\mt'$ is perturbed to a generic bundle $\mt''$with generic flags then $f(N)>1$ gives us at least two subbundles. One of these tends to $\mv$ as we degenerate $\mt''$ to $\mt'$ (along with the flags), and the others stay away from $\mv$.
This is because everything is transverse and smooth at $\mt'$ near $\mv$, so $\mv$ counts for $1$ point. The limit of the other point may only be a coherent subsheaf which we can saturate to a subbundle of $\mt'$ which will atleast as much  parabolic slope.
\end{remark}

In case the Grassmann dual is not polyrigid do the same argument but with
$\mv\oplus\mv\oplus\dots\oplus\mv\oplus\mq$ with $\mv$ repeated $N$ times.
\end{proof} 

\begin{remark}
It can be proved (corollary ~\ref{PolyRigidInv}) that polyrigidity is invariant under (Grassmann) duality. So the second part of the argument above is not needed.
\end{remark}

\begin{remark} A stable parabolic bundle is called rigid if any other
stable parabolic bundle with the same weights is isomorphic to it. Or the 
`corresponding' unitary representation\footnote{First use the transformation formulas to make the degree of the bundle $0$ and then norm the $\theta$'s so that they sum to $0$ at each point. To this modified data one can apply the theorem of Mehta and Seshadri.} of the fundamental group from the Mehta-Seshadri theorem is rigid.
\end{remark}
\begin{remark} This result points at a rigidity statement for essential inequalities in Eigenvalue problems. Namely what is the analogous rigidity statement for the essential inequalities in ~\cite{TW}?
\end{remark}
\begin{remark} The same argument shows the polyrigidity of any piece of the Harder Narasimhan filtration.
\end{remark}

There is a construction of intersection 1 situations due to Nori (private communication) and it will be useful to see if these situations correspond to polyrigid Witten bundles.
\begin{example}(M.V.Nori)
Let $V_1,\dots V_s$ be general subspaces of a vector space $W$ so that 
$$\sum_i dim(V_i)=dim(W)+1.$$

So, the natural map $\bigoplus V_i \to W$ has a $1$ dimensional kernel spanned
by $\sum v_i$ with $0\neq v_i\in V_i$. In Grassmannian $Gr(s-1,dim(W))$,
look for spaces $T$ which intersect all the $V_i$. There is exactly one such
space: $\sum_i \Bbb{C}v_i$. In terms of Schubert Calculus, let $I(a)$ be the subset of $\{1,\dots,dim(W)\}$ of cardinality $s-1$ given by $i(a)_1=rank(V_a)$ and $i(a)_l=dim(W)$ for $l\neq 1$. And the calculation above shows $\cup_a\sigma_{I(a)}$ to be the class of a point. One sees that the Witten bundle corresponding to this situation is rigid (and stable). Do all intersection 1 situations
arise from suitable generalizations of this theme?
\end{example}

\section{Moduli Spaces, Dimension Calculations}\label{irred}

If all our cycles are  \footnote{Recall that a cycle $I$ in $Gr(r,n)$ is said to be normalised if either $i_1>1$ or $i_r<n$.} normalised, then a situation with
$\langle\sigma_{I^1},\dots,\sigma_{I^{s}}\rangle_{d,D,n}\neq 0$
leads us to the moduli space of semistable degree $-d$ bundles, with parabolic  structure as in section ~\ref{HWitten}.  The moduli space classifies semistable pre-parabolic bundles (upto JH equivalence). Two parabolic bundles are isomorphic if there is a isomorphism which preserves parabolic structures, so the irrelevant parts of the flags are to be ignored (i.e we should pass to partial flags at each point
and ignore the pieces $F_i$ of the flags where the weights $\theta_i=\theta_{i+1}$ do not jump). So we get $s$ partial flags corresponding to each parabolic point. 

Let $\Lambda =\Lambda_1<\Lambda_2<\dots \Lambda_k= n$ be given and $V$ a vector space of dimension $n$. Let 
$$X_{\Lambda}(V)= \{F_{\Lambda_1}\subset F_{\Lambda_2}\dots \subset F_{\Lambda_k}=V\mid
dim(F_{\Lambda_l})=\Lambda_l\}$$

A simple dimension computation gives

$$dim(X_{\Lambda}(V))=\la_1(n-\la_1)+(\la_2-\la_1)(n-\la_2)+\dots+(\la_k-\la_{k-1})(n-\la_k).$$

For each choice of weights

$$\theta_1\geq\theta_2\geq\dots\geq\theta_n> \theta_1-1,$$
there is associated a partial flag variety. Namely pick the points
where this sequence jumps (as above). Let $\Lambda(\theta)$\footnote{$a\in\Lambda(\theta)\Leftrightarrow (a=n,\text{ or }\theta_{a+1}<\theta_a$).}  be the corresponding $\Lambda$. 

Now, let $\mt$ of degree $-D$ and parabolic weights $\theta^i_j$ be stable. We calculate the dimension of the corresponding moduli space $\mathcal{M}$ in terms of the partial  flag varieties associated to $\mt$.

it is clear what the dimension of the moduli space $\mathcal{M}$ should be 

$$dim(\text{Moduli of } \mt)=\sum _{l=1}^{l=s}dim(X_{\Lambda(\theta^l)})-dim(Aut(\mt))+1.$$

Over an open part of the moduli space the bundles stays isomorphic to $\mt$ because $\mt$ is a generic bundle. The product of the partial flag manifolds
 
$$X_{\Lambda(\theta^1)}(\mt_{p_1})\times X_{\Lambda(\theta^2)}(\mt_{p_2})\times\dots\times X_{\Lambda(\theta^s)}(\mt_{p_s})$$

gives a local universal family $\muu$ of parabolic bundles. Two such give the same point of the moduli space under the map $\muu \to \mathcal{M}$ if there is an element of $Aut(\mt)$ which takes one to the other. Hence the fiber over a given point of the moduli space is 

$$dim(Aut(\mt))-1$$

 because $\mt$ as a parabolic bundle has no automorphisms (is stable) other than scalars. For ease of notation let $\Lambda^l=\Lambda(\theta^l)$.The dimension of $X_{\Lambda}$ depends only on the multiplicities $\Lambda_1,\Lambda_2-\Lambda_1,\dots,\Lambda_k-\Lambda_{k-1}$. 

\begin{remark}There is a cyclic property of dimensions of partial flag varieties: let $\Lambda =\Lambda_1<\Lambda_2<\dots \Lambda_k= n$ be given, let
 $\Lambda' =\Lambda_1'<\Lambda_2'<\dots \Lambda_k'= n$ be given, with
$\Lambda_l'=\Lambda_{l+1}-\Lambda_{1}$ for $l<k$ and $\Lambda_k'=n$. Then, 
$dim(X_{\Lambda})= dim(X_{\Lambda'})$ from the formulas above. For instance let
$u_1,\dots u_{k-1}$ be the numbers $\Lambda_1,\Lambda_2-\Lambda_1,\dots, \Lambda_k-\Lambda_{k-1}$, then $dim(X_{\Lambda})=n^2-\sum_{i\leq j}u_iu_j$ and is invariant under the formula above. Closer to the point, these are dimensions of centralisers of elements in $U(n)$ with a central ratio.
\end{remark}

Finally we need the following simple result on automorphisms of generic bundles:
Let $\mathcal{V}=\mathcal{G}_{n,D}$. Clearly degree of $\mathcal{E}nd(\mathcal{V})=0$, so 

$$\chi(\mathcal{E}nd(\mathcal{V}))=n^2.$$ But if $\mathcal{V}$ is generic, all summands of
$\mathcal{E}nd(\mathcal{V})$ are either $\mathcal{O}(-1)$,$\mathcal{O}$ or $\mathcal{O}(1)$ all of which have $H^1=0$. So $dim(\mathcal{E}nd(\mathcal{V}))=n^2$. The identity is an automorphism and the the set of automorphisms is an open subset of the set of endomorphisms so is a group of dimension $n^2$. The following argument is needed to show that elements in $End(\mathcal{V})$ which give an automorphism
is a Zariski open set. It is not clear what groups are obtained this way.

\begin{lemma} Let X be a projective algebraic variety, $\mathcal{V}$ a vector bundle and
$x$ a point. Let $\mathcal{V}_x$ be the fiber of $\mathcal{V}$ at $x$.
We have a map $\phi_x:End(\mathcal{V}) \to End(\mathcal{V}_x)$.
\begin{enumerate} 
\item For $s \in End(\mathcal{V})$, the characteristic polynomial of $\phi_x(x)$ does not depend on $x$. 
\item Given  $s\in End(\mathcal{V})$, it is an automorphism if and only if  $\phi_x(s)$ is so.
\end{enumerate}
\end{lemma}
\begin{proof} The coefficients in characteristic polynomial are global regular
functions, hence constant. In particular the determinant is a constant.
\end{proof} 

We record our conclusion on the dimensions of these moduli spaces:
\begin{proposition}\label{CalculDim} Let $\mt$ be a stable parabolic bundle as above and assume that at each point the difference between the weights is $<1$ in absolute value. Let $X_{\Lambda_i}$ be the partial flag variety assigned to the parabolic point $p_i$ as above.

$$dim(\text{Moduli of } \mt)=\sum _{l=1}^{l=s}dim(X_{\Lambda_l})-n^2+1.$$
\end{proposition}

Assume now that the moduli of $\mt$ is $0$ dimensional and consists of a stable point (it is always connected). That means that there is a dense orbit of $Aut(\mt)$ acting on

$$X_{\Lambda_1}(\mt_{p_1})\times X_{\Lambda_2}(\mt_{p_2})\times\dots\times X_{\Lambda_s}(\mt_{p_s})$$

We will have(from the dimension of moduli of the parabolic bundle $\mt$)
$$0=\sum _{l=1}^{l=s}dim(X_{\Lambda(\theta^l)})-n^2+1.$$

On $\mt$ there is a remarkable other  choice of weights: namely to
$F_{\Lambda^l_k}$ assign the weight $rig^l_k=\frac{n-\Lambda^l_k}{n}$ these weights obviously lie in $[0,1)$. Denote the associated parabolic part of the weight by $\tilde{w}$  and parabolic slope by $\tmu$.

$$\tw(\mt)=\frac{1}{n}(\sum_{l=1}^{l=s}\sum_{a=1}^{k_l}(\Lambda^l_a-\Lambda^l_{a-1})(n-\Lambda^l_{a}))$$

(with the understanding that $\Lambda^l_0=0$)

which is

$$\frac{1}{n}(\sum _{l=1}^{l=s}dim(X_{\Lambda(\theta^l)}))= \frac{1}{n}(n^2-1).$$

so $\mu{\mt}=1-\frac{D}{n}-\frac{1}{n^2}$.

Suppose this structure is unstable, let $\ms$ be the Harder-Narasimhan element
which has $F_{\Lambda^l_a}\cap \ms$ is $u^l_a$ dimensional and of degree $-d$
and rank $r$. The uniqueness of the Harder-Narasimhan element gives
(the flags on $\mt$ are generic)

$$0=[r(n-r)+dn-Dr]-\sum_{l=1}^{l=s}\text{codimension}(Y_{l})$$

where

$$Y_{l}=\{S\in Gr(r,\mt_{p_l})\mid dim(F_{\Lambda^l_a}(\mt_{p_l}))\cap \ms_{p_l})\geq u^l_a; a=1,\dots,k_l\}.$$ 

($k_l$ is the length of $\Lambda^l$)
The codimension of $Y_l$ is 

$$r(n-r)-[u^l_1(\Lambda^l_1-u^l_1)+(u^l_2-u^l_1)(\Lambda^l_2-u^l_2)+\dots+(u^l_{k_l}-u^l_{k_l-1})(\Lambda^l_{k_l}-u^l_{k_l})$$

$=n(\tw_l(\ms))-nr+r(n-r)-\sum_{a=1}^{a=k_l}(u^l_a-u^l_{a-1})(-u^l_a)$

$=n(\tw_l(\ms))-r^2-\sum_{a=1}^{a=k_l}(u^l_a-u^l_{a-1})(-u^l_a)$
(with the usual understanding that $u^l_0=0$, and $\tw_l$ is the contribution to $\tw$ from $p_l$). So we get

$$-n(\tw(\ms))+sr^2+\sum_{l=1}^{l=s}\sum_{a=1}^{a=k_l}(u^l_a-u^l_{a-1})(-u^l_a)+dn-Dr+r(n-r)=0$$

or that 
$$(\tw(\ms)-d)n=sr^2+r(n-r)+\sum_{l=1}^{l=s}\sum_{a=1}^{a=k_l}(u^l_a-u^l_{a-1})(-u^l_a)-Dr.$$

We have 
$\tw(\ms)-d > \tmu(\mt) r$ by assumption, so

$$sr^2+r(n-r)+\sum_{l=1}^{l=s}\sum_{a=1}^{a=k_l}(u^l_a-u^l_{a-1})(-u^l_a)-Dr=(\tw(\ms)-d)n>rn-Dr-\frac{rd}{n}$$

So
$$(s-1)r^2+\sum_{l=1}^{l=s}\sum_{a=1}^{a=k_l}(u^l_a-u^l_{a-1})(-u^l_a)>-\frac{r}{n}.$$

But we claim
\begin{equation}\label{RIG2}
(s-1)r^2+\sum_{l=1}^{l=s}\sum_{a=1}^{a=k_l}(u^l_a-u^l_{a-1})(-u^l_a)+1\leq 0.
\end{equation}
because `$\ms$ is rigid' (as will be explained below) being the Harder-Narasimhan element in a rigid $\mt$. And hence a contradiction.

First note that $\ms$ gets a partial flag in each of its  fibers $\ms_{p_l}$, the $l$th one being made out of $F_{\Lambda^l_a}(\mt_{p_l}))\cap \ms_{p_l}$
let this corresponding flag be in $X_{\Delta_l}(S_{p_l})$
where $\Delta_l= \Delta^l_1<\dots<\Delta^l_{b_l}$ is the set $\{u^l_a\mid a=1,\dots k_l\}$

For fix generic bundles $\mt$, and generic (partial flags) in
$$fl_0=(F^1,\dots,F^s)\in FL(\mt)=X_{\Lambda_1}(\mt_{p_1})\times X_{\Lambda_2}(\mt_{p_2})\times\dots\times X_{\Lambda_s}(\mt_{p_s})$$ and a $\ms\subset\mt$ in

$$\cap_l Y_l(F^l)$$

(there is a unique one by assumption, because it a Harder-Narasimhan element for a choice of weights by assumption) where

$$Y_{l}=\{\ms\in Gr(r,\mt_{p_l})\mid dim(F_{\Lambda^l_a}(\mt_{p_l}))\cap \ms_{p_l})\geq u^l_a; a=1,\dots,k_l\}.$$ 

($k_l$ is the length of $\Lambda^l$)

Now look at 

$$\mathcal{J}=\{(G^1,\dots,G^s)\in FL(\mt)\mid \ms \in \cap_l Y_l(G^l)\}.$$

We claim there is a map 

$$\phi:\mathcal{J}\to FL(\ms)$$

where

$$FL(\ms)=X_{\Delta_1}(\ms_{p_1})\times X_{\Delta_2}(\ms_{p_2})\times\dots\times X_{\Delta_s}(\ms_{p_s}).$$

The map is the obvious one given flags on $\mt_{p_l}$, intersect them with
$\ms_{p_l}$. The map is surjective near $fl_0$ (the original flags on $\mt$).
Now Aut(T) acts transitively on the elements of $FL(\mt)$ (at least the ones near $\mt$) because the original choice of weights $\theta^l_a$ gave a stable rigid
parabolic bundle. So given $fl_1,fl_2\in FL(\mt)$, there is an element
in $A\in Aut(\mt)$ with $A(fl_1)=fl_2$. Now $\ms$ is in  

$$\cap_l Y_l(G^l)$$

and in

$$\cap_l Y_l(H^l)$$

where $fl_1=(G^1,\dots,G^s)$ and $fl_2=(H^1,\dots,H^s)$ and is the unique
such element. Since $u(fl_1)=fl_2$, we have $u(\ms)=\ms$. Putting it all together we get that $Aut(\ms)$ acts transitively on $FL(\ms)$ and the equation in the lemma is in contradiction with equation ~\ref{RIG2}.

\begin{lemma} Let $\ms$ be a generic bundle of degree $-d$, rank $r$. Let

$\Delta_l= \Delta^l_1<\dots<\Delta^l_{b_l}=r$ for $l=1,\dots,s$. Let

$$FL(\ms)=X_{\Delta_1}(\ms_{p_1})\times X_{\Delta_2}(\ms_{p_2})\times\dots\times X_{\Delta_s}(\ms_{p_s}).$$

Suppose $p\in FL(\ms)$ and $Aut(\ms)$ acts transitively on $FL(\ms)$

Then 

$dim(Aut(\ms))-1 \geq \sum dim(X_{\Delta_l})(\ms_{p_l}).$

\end{lemma}

\begin{proof} If a group acts transitively on a variety $X$ (or has a dense open orbit), then the dimension of the group is at least as much as that of $X$.
Here the group is $Aut(\ms)/\Bbb{C}^*$ and the variety $X=FL(\ms)$.
\end{proof}
Hence we have shown that the new structure on $\mt$ is semistable. To summarise
the results of this section.

\begin{theorem}Let $\mt$ be a generic bundle of degree $-D$ and rank $n$
endowed with $s$ complete (and generic) flags and weights $\theta^l_a: l=1,\dots s; a=1,\dots n$ so that for any $l$,

$$\theta^l_1\geq\theta^l_2\geq\dots\geq\theta^l_n >\theta^l_1-1.$$(No other assumptions). Suppose $\mt$ with this structure is stable. We have

$$dim(\text{Moduli of } \mt)=\sum _{l=1}^{l=s}dim(X_{\Lambda(\theta^l)})-n^2+1.$$

Suppose this dimension is $0$. Then define a new set of weights by assigning to
$F_{\Lambda^l_k}(\mt_{p_l})$ assign the weight $rig^l_k=\frac{n-\Lambda^l_k}{n}$.Then $\mt$ with this new weights is semistable of total slope $1-\frac{D}{n}-\frac{1}{n^2}$. Call this new parabolic bundle $Rig(\mt)$ (which can be defined even if $\mt$ is not stable).
\end{theorem}
\begin{remark} Note the following additive property of $Rig(\mt)$: 

$$Rig(\oplus_{i=1}^N \mv)=\oplus_{i=1}^N Rig(\mv).$$
\end{remark}
We remark that this new structure is stable as well. For the integrality conditions force stability=semistability. If $\ms$ is a contradiction ($r=rank(\ms)$) to stability
we will have ($\mu$ is the rigidity weight)

$\mu(\mt)=1-\frac{D}{n}-\frac{1}{n^2}$

It is clear that (the weights have denominators $n$), that

$rn\mu(\ms)\in \Bbb{Z}$. But

$$rn\mu(\ms)=nr -Dr-\frac{r}{n}$$

and the last element is clearly not an integer. For later reference
note the result
$rn\mu(\ms)< nr -Dr-\frac{r}{n}$,
so $rn\mu(\ms)\leq nr -Dr-1.$
Hence
\begin{equation}\label{STRONG}
\mu(\ms)\leq 1 -\frac{D}{n}-\frac{1}{rn}
\end{equation}
\begin{example}\label{egg}
In this example we look in $Gr(8,12)$, $d=D=0$, and cycles
\begin{enumerate}
\item $I^1=\{3,4,5,7,8,10,11,12\}.$
\item $I^2=\{2,3,5,6,8,9,11,12\}.$
\item $I^3=\{2,3,5,6,8,9,11,12\}.$
\end{enumerate}
Then the Witten Moduli Space $\mathcal{M}$(which is generically stable)  is $6$ dimensional and the intersection number (which equals $H^0(\mathcal{M},\Theta)$) $ = 6$). So $\Theta$ (which is ample) is not base point free and does not give an embedding generically. The Grassmann dual Witten Bundle is stable in this example. The example was found (and checked) by a computer program.  
\end{example}
\begin{example}\label{EXT} To illustrate the importance of extensions of simple factors in parabolic moduli spaces, consider the Witten bundle of the intersection $d=D=0$, in $Gr(5,9)$.
\begin{enumerate}
\item $I^1=\{3,5,6,8,9\}.$
\item $I^2=\{2,4,5,7,9\}.$
\item $I^3=\{2,3,5,8,9\}.$
\end{enumerate}

In this case the parabolic moduli space is a point, but the product of the $3$ partial flag varieties is not acted upon transitively (generically) by $Sl(5)$. the flag varieties have dimensions $4+4$, $4+4+1$ and $6+2$ respectively and the dimension of the group is $24$.  One checks that the line $\ml$ in this five dimensional space with Schubert position $\{5\}$, $\{4,\}$ and $\{2,\}$ contradicts stability. In fact we can show that $\mv$ is an extension of rank 1 stable bundles (which shows that the parabolic moduli space is a point).
\end{example}

\section{Irredundancy}\label{I1}
Consider a polyrigid inequality $Ineq(r,d,D,I^1,\dots,I^s,n)$. Choose generic flags on $\mg_{D,n}$ and let $\mv\subset\mg_{D,n}$ be the point corresponding to 
$$\langle\sigma_{I^1},\dots,\sigma_{I^{s}}\rangle_{d,D,n}=1$$

We will assume to start with that the Grassmann dual Witten bundle of $\mv$ is normalised. This implies that either $i^a_1=1$ or $i^a_r=n$ for each $a\in\{1,\dots,s\}$.

The main induction step is the following theorem which we assume in this section and prove it in section ~\ref{JHSaga}:

\begin{theorem}\label{UneTheoreme}
We can find $\mw\subset \mt=\mg_{D,n}$ so that
\begin{enumerate}
\item $\mv\subset\mw$.
\item $dim(\mw,\mt)=0$ and is a polyrigid intersection.

\item The Grassmann dual Witten bundle of $\mw$ is the direct sum of the same simple factor.

\end{enumerate}

It then follows that,
\begin{enumerate}
\item The Grassmann dual bundle of $\mw$ is normalised. This follows from the assumption that the Grassmann dual Witten bundle of $\mv$ is normalised.
\item $\mv$ is a contradiction to stability (not to semistability) of the
Witten Bundle of $\mw$ (this is automatic from the previous step).
\end{enumerate}
\end{theorem}

Assuming this theorem we can immediately show the strong irredundancy of polyrigid inequalities:

\begin{defi}\label{min} Let $\mt$ be a vector bundle on $\Bbb{P}^1$ with complete flags on its fibers at $p_1,\dots,p_s$. Let $\mv\subset\mt$ be a subbundle and it recieves complete flags on its fibers. Suppose we have a pre-parabolic
structure on $\mv$, then there is a minimal extension of the parabolic structure to all of $\mt$. To define this we take a point $p=p_k$, and suppose 
$\mv_p\in \Omega^o_I(F_{\bull})(\mt_{p})$, where $I=\{i_1,\dots,i_r\}$ and the weight assigned to $F_l(\mv_p)$ is $\theta_l$. Then if ($r=rank(\mv),n=rank(\mt)$).
\begin{enumerate}
\item $i_u<j< i_{u+1}$ assign weight  $\theta_{u+1}$ to $F_{j}(\mt_p)$.
\item $0<j< i_1$ then assign the weight $\theta_{1}$ to $F_{j}(\mt_p)$.
\item $i_r<j\leq n$ assign the weight $\theta_1-1$ to $F_{j}(\mt_p)$.
\item $j=i_u$ then assign the weight $\theta_{u}$ to $F_{j}(\mt_p)$.
\end{enumerate}
The induced parabolic structure on $\mv$ is the old one.

\end{defi}

We want to consider the minimal extension of the Witten parabolic structure on $\mw\subset\mt$. The first claim is that the induced weights on $\mq=\mt/\mw$ can also be obtained as follows:
\begin{enumerate}
\item Form the Witten bundle of $\mq^*\subset\mt^*$ ($\mt^*$ has natural induced flags on it). Take the dual of this structure and get one on $\mq$.
\item Form the Rigidity bundle of this  parabolic structure on $\mq$.
\end{enumerate}

\begin{claim}The Rigidity structure on $\mq$ in (2) above is same as the weights induced on $\mq$ by the minimal extension of the Witten structure on $\mw$ to $\mt$.
\end{claim}
\begin{proof}
Let $W\subset T=\Bbb{C}^n$, be in Schubert position $I=\{i_1<i_2<\dots<i_r\}$
with respect to the complete flag $F_{\bull}$ on $T=\Bbb{C}^n$. We get the induced flag $F(W)_{\bull}$ on $V$ and the flag $F_{\bull}(Q)$ on $Q=T/W$. There are two ``Horn'' type situations: $W\subset T$ and $Q^{*}\subset T^{*}$. There are resulting weights on $W$, $Q^{*}$ and hence $Q$ from this situation. 

Let $J=dual(I,n)=\{a\in\{1,\dots,n\mid n+1-a \not\in I\}$, and let
$J=\{j_1<\dots<j_{n-r}\}$. We are assuming that $Q^*$ has a normalised Schubert position so either $j_1 \neq 1$ or $j_r \neq n$.
\begin{enumerate}

\item Weights on $W$: On $F(W)_a$ assign the weight $1+\frac{1}{n-r}(a-i_a)$. This structure is what one gets from the Horn conjecture (the Witten bundle). 

\item Weights on $Q^*$: on $F(Q^*)_a$ assign the weight $1+\frac{1}{r}(a-j_a)$. 

\item \label{RIG}Weights on $Q$: on $F(Q)_a$ assign the weight $1-\frac{1}{r}((n+1-a)-j_{n+1-a})$. by  dualising the previous one. The weight of $F(Q)_{a}$ is different from weight of $F(Q)_{a+1}$ if and
only if  

$$(n+1-a)-j_{n+1-a}\neq n-a-j_{n-a}$$

or that

$$j_{n+1-a}-j_{n-a}\neq 1$$

Let the complement of $I$ in $\{1,\dots,n\}$ be $\{\alpha_1<\dots<\alpha_{n-r}\}$. In this notation $j_{n+1-a}=n+1-\alpha_a$. So 
$\alpha_{a+1}>1+\alpha_a$. If $\alpha_a+1=i_b$, then $a=i_b-b$.

\item Weights on $T$: to the $F_{i_a}$ assign weight $1+\frac{1}{n-r}(a-i_a)$, to all others assign the smallest weight posible. That is if $i_{a-1}< l \leq i_{a}$ assign the weight $\frac{1}{n-r}(a-i_a)$ to $F_l(T)$. If $i_r <n $ (in which case $i_1=1$, because the Grassmann dual is normalised) assign the weight $0$ to $F_n$.

 This choice of weights
induces weights on $Q$. The induced weight on $F(Q)_{i_a-a}$ is $1+\frac{a-i_a}{n-r} $. This last weight is the rigidity weight from ~\ref{RIG} on $Q$. If $i_r<n$, the weight assigned from the above procedure on $F_{n-r}(Q)$ is $0$ which agrees with the rigidity weights.
\end{enumerate}
\end{proof}
\begin{example} We give an example on how weights are assigned, so as to also
clarify the freedom we have in the extension of the weights from $\mv$ to $\mt$. We do the point case: Let $\mv$ be a $4$ dimensional space in $\Bbb{C}^8$. Let the Schubert position of $\mv$ be $\{2,4,6,8\}$ with respect to the flag $F_{\bull}(\Bbb{C}^8)$.
So to $F_2(\Bbb{C}^8)$,$F_4(\Bbb{C}^8)$,$F_6(\Bbb{C}^8)$,$F_8(\Bbb{C}^8)$, we assign the weights $\frac{8-4+1-2}{4}=\frac{3}{4}$, $\frac{2}{4},\frac{1}{4}$ and $0$ respectively. To $F_7$ the conservative extension is $0$, but we can go as far as $\frac{1}{4}$ and still have a valid choice of weights.
\end{example}
\begin{remark}\label{Thestretch}The freedom in the above construction: Although we chose the most conservative extension so as to minimise the slope of the quotient, we will need to know the maximum slope the quotient can acquire so that the subbundle has the  weights as in the Horn Conjecture. An easy calculation shows that we can add $\frac{1}{rank(\mq)}$ to each of the weights not appearing in the subbundle at any parabolic point. So that process can lead to increase of slope of $\frac{1}{rank(\mq)}$ for each parabolic point. If the Witten bundle of $\mq^*\subset \mt^*$ is $M$ copies of the same (parabolic) bundle $\ml$ - then the freedom in the above construction is $\frac{1}{rank(\ml)}$. 
\end{remark}

Now (the global case) take the weights obtained by Horn Conjecture on $\mw$ and take the most conservative extension $\mu$ to $\mt$. The
induced weights on $\mq$ come from the rigidity weights coming from the structure on $\mq$ obtained by taking the dual of the Horn conjecture weights on $\mq^{*}$. Let $\mu$ denote the slopes in this conservative extension.

Let $\mq=\mt/\mw$ which as a Grassmann dual Witten bundle is a direct sum of $M$(say) copies of the same bundle $\ml$ of rank $l$. So its Rigidity bundle has slope $1+degree(\mq)/rank(\mq)
-\frac{1}{l^2}$. We have the following formulas for the slopes $\mu$.
\begin{enumerate}
\item $\mu(\mw)=\mu(\mv)=1+degree(\mq)/rank(\mq)$.
\item $\mu(\mq)=degree(\mq)/rank(\mq)-\frac{1}{l^2}$.
\end{enumerate}

Therefore we have that $\mv$ contradicts semistability (and hence satisfies condition (1) of lemma ~\ref{lemmacondition}). We have the freedom of adding a constant $c<\frac{1}{l}$ to all the weights of $\mq$ at any parabolic point. We add the same constant $c$ to all the weights of $\mq$ (at every parabolic point) so that the slope of $\mq$ becomes equal to that of $\mw$. The constant $c$ is given by $sc=\frac{1}{l^2}$, or $c=\frac{1}{sl^2}$. We have a structure on $\mt$ for which the slope satisfies
$\mu'(\mw)=\mu'(\mv)=\mu'(\mq)=1+degree(\mq)/rank(\mq)$.

  At this choice of weights, the weights of $\mw$ and $\mq$ are separated. Now perturb \footnote{This is a  small perturbation with differences $<< \frac{1}{sl^2}$} the weights on $\mw$ so that it remains semistable, the  subbundle $\mv$  still contradicts stability and the weights of $\mv$ and $\mw/\mv$ are now separated out. We have used induction in the last step.  Then, add (or subtract) a small constant to the weights in $\mq$ so that $\mu''(\mw)=\mu''(\mq)$. So we have satisfied condition (2) in lemma ~\ref{lemmacondition}. This finishes the last step.

\section{Examples}\label{examples}
All the examples below are from the classical part. One could use the transformation formulas to get ``quantum'' examples - but we expect that there are quantum polyrigid situations that do not arise from the classical ones by application
of the transformation formulas (to do that we need to look through quantum examples which we havent done, need some technology to code the quantum Littlewood Richardson rule). 

In each of these examples we are given an intersection in a $Gr(r,n)$ corresponding to an intersection of Schubert cells indexed by $I$'s. Having chosen generic flags in $\mt=\Bbb{C}^n$, we will find a subspace $\mv$ of rank $r$ in the Schuberts position corresponding to the $I$'s. Let $\mq=\mt/\mv$ and $\mq^*\subset \mt^{*}$ in the Grassmann dual Schubert position with respect to the induced flags on $\mt^*$ - denoted by $J$'s.

\begin{example}\label{example1} For $n=2,r=1,d=D=0,s=3$ there is really one inequality which needs to be shown to be irredundant: $I^1=\{1\}$, $I^2=I^3=\{2\}$. The Grassmann dual Horn bundle is polyrigid in this case. The Horn weights are on a one dimensional trivial vector bundle $\mv$ on $\Bbb{P}^1$. The Horn weight at $p_1$ is $1$,at $p_2,p_3$ is $0$. So the conservative extension of these weights to $\mt=\mg_{0,2}$ is $(0,0)$ at $p_2,p_3$  and $(1,0)$ at $p_1$. The only contradiction to semistability with this choice of weights is $\mv$.
\end{example}

\begin{example}\label{example2} For $n=8, r=5, d=D=0, I^1=\{3,4,5,7,8\},I^2=I^3=\{2,3,5,6,8\}$. The Grassmann dual is happenning in $Gr(3,8)$, with $J^1=\{3,7,8\}$,\\$J^2=J^3=\{2,5,8\}$. The Grassmann dual Witten Bundle bundle is stable and rigid. The Witten weights on the bundle $\mv$ coming from $Gr(5,8)$ are $(\frac{1}{3},\frac{1}{3},\frac{1}{3},0,0)$, and $(\frac{2}{3},\frac{2}{3},\frac{1}{3},\frac{1}{3},0)$ at $p_2$ and $p_3$. When induced up to $\mt=\mg_{0,8}$ (conservative extension), the weights are $(\frac{1}{3},\frac{1}{3},\frac{1}{3},\frac{1}{3},\frac{1}{3},0,0,0)$, and $(\frac{2}{3},\frac{2}{3},\frac{2}{3},\frac{1}{3},\frac{1}{3},\frac{1}{3},0,0)$ at $p_2$ and $p_3$. And the quotient $\mq=\mt/\mv$ has (induced)weights -$(\frac{1}{3},\frac{1}{3},0),(\frac{2}{3},\frac{1}{3},0),(\frac{2}{3},\frac{1}{3},0)$ of slope $1-\frac{1}{9}$.

There is one more choice of weights on $\mq$. The dual $\mq^*\subset\mt^*$ of $\mq$ has Witten weights: $(\frac{3}{5},0,0), (\frac{4}{5},\frac{2}{5},0),$ and $(\frac{4}{5},\frac{2}{5},0)$. This induces weights on $\mq$ (we take the negatives and $1$ and rearrange in ascending order)- $(1,1,\frac{2}{5}),(1,\frac{3}{5},\frac{1}{5})$ and $(1,\frac{3}{5},\frac{1}{5})$. The associated rigidity weights of this structure on $\mq$ is $(\frac{1}{3},\frac{1}{3},0)$, $(\frac{2}{3},\frac{1}{3},0)$ and $(\frac{2}{3},\frac{1}{3},0)$. This is the same as the structure at the end of last paragraph.

From the theory developed in the paper, $\mv$ is the Harder-Narasimhan maximal contradictor of stability for $\mt=\mg_{0,8}$. But there are other contradictors of stability because $\mv$ is not stable. Example the line $\mathcal{S}\subset\mv$ in Schubert position (with respect to the induced flags on $\mv$)-$\{5\},\{4\},\{2\}$. This picks up weights: $0$,$\frac{1}{3},\frac{2}{3}$, and hence has the same slope as that of $\mv$ ($=1$).

So how does the paper make $\mv$ the only contradictor of semistability? First we add weights $\epsilon $ to all the weights of the quotient so that its slope becomes equal to that of $\mv$. In this case, the constant addition $\epsilon$ is given by $9\epsilon =\frac{1}{3}$. So now the weights of $\mt$ become $(\frac{1}{3}+\epsilon ,\frac{1}{3}+\epsilon ,\frac{1}{3},\frac{1}{3},\frac{1}{3},\epsilon ,0,0)$, and $(\frac{2}{3}+\epsilon ,\frac{2}{3},\frac{2}{3},\frac{1}{3}+\epsilon ,\frac{1}{3},\frac{1}{3},\epsilon ,0)$ at $p_2$ and $p_3$.  One checks in general that there is ``space'' for this addition (so that the order of the weights remains unchanged). The weights of $\mv$ and of $\mt/\mv$ are now separated (the absolute value of the difference of any weight from $\mv$ and one from $\mt/\mv$ is greater than $0$ and less than $1$). This new structure on $\mt$ is semistable. If we can perturb it so that $\mv$ and $\mq$ become stable and of the same slope and so that the weights are now in the interior of the large weightspace , then we are done. This is because then,
\begin{enumerate}
\item $\mv$ and the quotient (if $\mv$ splits off) are the only contradictors
of semi-stability.
\item we can join this set of weights to the conservative extension of the Horn weights on $\mv$ and head out just a little. For a short time the subs which were not contradicting semistabilty continue to not do so, and $\mv$ immediately starts to contradict stability.
\end{enumerate}

Such a perturbation is possible by induction. We can perturb the weights on $\mv$ and on $\mq$ so that they become stable and the weights are now in the interior. But we have changed the slopes of each (we cannot perturb keeping the slopes the same and remain in the large weight space, because the initial point we have is not in the interior). But the separation of the weights of $\mv$ and of $\mq$, allows us to correct this by adding a small constant to the weights of $\mq$ for example.
\end{example}

\begin{example}\label{example3} We give an example for the  Grassmann dual Witten bundle  not stable (but polyrigid).
$Gr(2,5), I^1=I^2=I^3=\{2,5\}$. The dual Grasmann intersection is in $Gr(3,5)$
and corresponds to $(J^1=J^2=J^3=2,3,5)$. The Witten bundle for $\mq^*$ is unstable - there is a line which contradicts semistability in Schubert position $\{3,\},\{2\},\{2\}$. This gives a subspace of dimension $2$ of $\mq$ in Schubert position $\{2,3\},\{1,3\},\{1,3\}$. This in turn gives a $4$ dimensional subspace of $\mt$ in Schubert position $\{2,3,4,5\},\{1,2,4,5\},\{1,2,4,5\}$.. This hyperplane contains $\mv$ and is unique in its Schubert class. Now consider the Horn bundle of $\mw$: weights $(0,0,0,0),(1,1,0,0),(1,1,0,0)$. The conservative extension of thse weights to $\mt$ is $(0,0,0,0,0),(1,1,0,0,0), (1,1,0,0,0)$. One checks that $\mw$ is the Harder-Narasimhan bundle of $\mt$, and $\mv$ has the same slope of $\mw$ ($=1$). There are two issues here: The first to make $\mv$ the Harder-Narasimhan bundle, and also to produce weights for which $\mv$ s the sole contradiction to semistability. 

We can first add a constant to the weights of $\mt/\mw$ so that $\mt$ becomes semistable: the constant is $3\epsilon =1$, and the weights on $\mt$ are now

$$(\epsilon,0,0,0), (1,1,\epsilon,0,0),(1,1,\epsilon,0,0).$$

With these weights $\mt$ is semistable and $\mw$, $\mv$ (among others) contradict stability. The weights on $\mw$ and on $\mt/\mw$ are separated.

So we perturb the weights on $\mw$ so that $\mv$ is the sole contradictor to stability (and the weights are in the interior of the large weight space) inside $\mw$. The separation of the weights (of $\mw$ and $\mt/\mw$) says that we can then adjust the weights on $\mt/\mw$ so that the slope of $\mt/\mw$ and that of $\mw$ become the same. We now appeal to lemma ~\ref{lemmacondition}.

\end{example}

\section{Proof of Theorem ~\ref{UneTheoreme}}\label{JHSaga}

Using the notation of Theorem ~\ref{UneTheoreme}, we find that we have nothing to prove if the Grassmann dual Witten bundle of $\mv$ is stable or is direct sum of the same simple parabolic bundle. We have to take care of the case when the Grassmann dual Witten bundle of $\mv$ is not of this type. We will assume to start with that the Grassmann dual Witten bundle of $\mv$ is normalised. This implies that either $i^a_1=1$ or $i^a_r=n$ for each $a\in\{1,\dots,s\}$.

The following theorem shows that if $\mv$ is a polyrigid semistable bundle with only one isomorphism type of stable parabolic bundle in its Jordan-H\"{o}lder series, then it is a direct sum of (several copies) of this stable bundle.

\begin{lemma}\label{Rigid3}
Let $\mv$ be a polyrigid stable bundle \footnote{That is $\mv$ is  rigid and hence has generic underlying bundle.} (with pairwise difference of weights at any point $<1$). Then, $\mv\oplus\dots\oplus\mv$ ($N$ times)has no deformations as a (semistable) parabolic bundle.

\end{lemma}
\begin{proof} Look at a $\mv\subset\mw=\mv\oplus\dots\oplus\mv$ as a coordinate factor. One easily checks that the associated point on the relevant Quot scheme is smooth. This follows from the genericity of the underlying bundle of $\mv$. 

We just need to check that this subbundle deforms with any deformation of $\mw_0$ (preserving its Schubert position). That is the tangent space of subbundles of the same (Schubert) type as $\mv$ of $\mw_0$ at $\mv$ is the expected dimension. This follows from the rigidity condition on $\mv$. We also have that
the subbundle $\mv$ does not deform as a parabolic bundle (because it is rigid). So combining these we have $N$ factors of $\mv$ in any nearby $\mw_t$. The direct sum map is an isomorphism at $\mw_0$, hence for nearby $t$. So we just need to show that the tangent space of subs of the same Schubert position as $\mv$ is of the expected dimension. 

The actual dimension of the tangent space is (the HOM is homomorphism of parabolic bundles)

$$HOM(\mv,\oplus_{i=2}^{N}\mv),$$

which is clearly $N-1$ dimensional. The expected dimension is to be calculated:

Point calculation: Let $i^l_1<\dots <i^l_{k_l}=r$ be the associated partial flag variety of $\mv$. The $l$th Schubert cell of $\mv$ in $\mw_0$ is

$$\{V\in \Bbb{C}^{Nr}\mid dim(V\cap F_{Ni^l_a}(\Bbb{C}^{Nr}))=i^l_a\}$$

which has dimension

$$i^l_1(Ni^l_1-i^l_1)+ (i^l_2-i^l_1)(Ni^l_2-i^l_2)+\dots+(i^l_{k_l}-i^l_{k_{l-1}})(Ni^l_{k_l}-i^l_{k_l})$$

=$$(N-1)[i^l_1 i^l_1+ (i^l_2-i^l_1)i^l_2+\dots+(i^l_{k_l}-i^l_{k_{l-1}})i^l_{k_l}].$$

So {\em {expected}} dimension is 

$$(N-1)[-(s-1)r^2 +\sum_l(N-1)[i^l_1 i^l_1+ (i^l_2-i^l_1)i^l_2+\dots+(i^l_{k_l}-i^l_{k_{l-1}})i^l_{k_l}]=$$

$$(N-1)[1-dim(\text{moduli of }\mv)],$$

from the proposition ~\ref{CalculDim} and which is $(N-1)$ in our situation.
\end{proof}
\begin{remark}\label{ALTERNATE} Parabolic bundles $\mv$ on a curve $C$ can be thought of as $\Gamma=\frac{\Bbb{Z}}{N\Bbb{Z}}$ equivariant vector bundles $\mw$ on a Galois cover $X$ of $C$ with Galois group $\Gamma$. The deformation space is therefore expected to be $(H^1(X,End(\mw))^{\Gamma}$. It is clear with this formalism that if the deformation space of $\mw$ is $0$ dimensional then there are no deformations of $\oplus\mw$ as well. But one has to develop enough homological  algebra and deformation theory to carry out this argument. 
\end{remark}

\begin{remark} Let $\ml$ be a rigid local system on $\Bbb{P}^1-S$ for $S$ a finite set in the sense of N. Katz. If $\ml$ is unitary, then the above shows that
$\oplus_{i=1}^N \ml$ has no deformations as a unitary representation and hence in the space of all representations (because of  Weil's tangent space calculations ~\cite{weil}). The same result is true for any rigid local system $\ml$, again because of  Weil's results. The result shown is slightly stronger because, parabolic bundles and unitary representations coincide only up to extensions.
\end{remark}

We prove a result which allows us to get inductive grip on semistable but not stable bundles. The theorem is the following:
\begin{theorem}\label{JHH} Let $\mf$ be a semistable but not stable bundle whose weights at any point have differences less than $1$ in absolute value. There exists a
 subbundle  $\mg$ (resp a quotient $\mq$) of $\mf$ which satisfies
\begin{enumerate}
\item The slope(induced) of $\mg$ (resp $\mq$) is the same as that of $\mf$.
\item It has no infinitesimal deformations inside of $\mf$(a deformation which does not change its Schubert position with respect to the partial flag varieties corresponding to the parabolic fibers of $\mf$).
\item The Jordan-H\"{o}lder Filtration of $\mg$ (resp. of $\mq$)has exactly one simple factor.
\end{enumerate}
\end{theorem}
\begin{proof} The proof is quite general and applies also to moduli of vector bundles on a curve. In the parabolic case we ignore parts of the flags on parabolic bundles where the weights do not jump, and assume that the absolute values of differences of weights at any parabolic point are $<1$. We prove only the statement for subbundles. The case for quotients follows by duality.

Let $\mathcal{A}$ be one of the following
\begin{enumerate}
\item The abelian category of slope semistable bundles of a given slope $\mu$ on a smooth projective curve.
\item The abelian category of slope semistable parabolic bundles of a given
parabolic slope over a fixed curve and fixed parabolic points.
\end{enumerate}

In this situation by classical Jordan-H\"{o}lder theory, for any $\mathcal{F}\in \mathcal{A}$, there is a filtration

$$\mathcal{F}=\mathcal{F}_0\supset\mathcal{F}_1\supset\dots\supset\mathcal{F}_s\supset\mathcal{F}_{s+1}=0$$

with the graded quotients

$$\frac{\mathcal{F}_l}{\mathcal{F}_{l+1}}$$

simple (i.e stable) for $l=0,\dots,s$. Such a filtration is not unique but the graded quotients are. But some uniqueness can be obtained as follows. Let
$\mathcal{F}\in \mathcal{A}$ and let $\mathcal{V}$ be a simple object in $\mathcal{A}$.
Consider a maximal element $\mathcal{G}$ of the set

$$\{\mathcal{G}\subset \mathcal{F}\mid \text{ each graded quotient in the JH filtration of } \mathcal{G}\text { is isomorphic to } \mathcal{V}\}$$ 

Claim: $\mathcal{G}$ has no infinitesimal deformations in $\mathcal{F}$ (in the parabolic case the infinitesimal deformation should be one preserving the parabolic structure).

If it did we will find a map $\mathcal{G}\to \frac{\mathcal{F}}{\mathcal{G}}$
in the category (tangent spaces to Quot-Schemes even in the 'obstructed case').
The image is a $\mathcal{L}\subset\frac{\mathcal{F}}{\mathcal{G}}$ all whose
Jordan-H\"{o}lder graded pieces are isomorphic to $\mathcal{V}$ (being a surjective image of $\mathcal{G}$). Let $\tilde{\mathcal{L}}$ be the inverse-image in $\mathcal{F}$, this contains $\mathcal{G}$ and all the graded pieces in the JH filtration are isomorphic to $\mathcal{V}$, contradiction. 
\end{proof}

\begin{remark} Let $\mathcal{V}$ be a parabolic bundle such that the associated moduli space is $\mathcal{M}$. We will call $\mv$ polyrigid if $\mathcal{M}$ is a point. But as a vector bundle with flags it might still admit deformations.
For example one could vary along the Ext groups (see example ~\ref{EXT}). We will say that $\mv$ has no deformations as a parabolic bundle if it has no deformations as a vector bundle with flags at the parabolic points (and where as usual we ignore parts of the flags where the weights do not jump). 
\end{remark}

In the case $\mf$ is polyrigid, we have stronger properties of $\mg$ produced above.

\begin{lemma}\label{StrongRigid}
Now suppose $\mathcal{F}$ is a polyrigid semistable (with generic flags and generic underlying bundle). Let $\mathcal{G}$ be as above and assume $\mathcal{G}\neq 0$. We claim that the Witten bundle of $\mathcal{G}\subset\mathcal{F}$ is polyrigid (in fact has no deformations as a parabolic bundle).
\end{lemma}
\begin{proof} First note that under the conditions the parabolic bundle $\mg$
although not stable has no deformations as a parabolic bundle. This follows from
lemma ~\ref{Rigid3} above. The Witten Bundle for $\mg\subset \mathcal{F}$ has
associated flag varieties which are obtaining by (possibly) ignoring some parts of the partial flag varieties associated to $\mg$ as a parabolic bundle. Hence the Witten bundles does not have any deformations as a parabolic bundle. Let us be explicit about the partial flag varieties:

Point calculation: If the partial flag variety of $\mt$ at $p_l$ is ($rank(\mt)=n)$, $i^l_1<\dots<i^l_{k_l}=n$ and 
suppose $dim(\mg_{p_l}\cap F_{i^l_a}(\mt_{p_l}))=u^l_a,$

The partial flag variety for $\mg$ as a parabolic bundle (at $p_l$) is 

$$F_{\lambda_1}\subset F_{\lambda_2}\subset\dots\subset F_{\lambda_u}=\mg_{p_l}$$

where the $\lambda$'s are in the set $\{u^l_a\}$. To $F(\mg)_{u^l_a}$ the weight assigned is $\frac{n-r+u^l_a-i^l_a}{n-r}$ if $u^l_a\neq u^l_{a-1}$, so the jumps can occur only in the set $\{u^l_a\}$.

\end{proof}

We next prove that direct sum of several copies of a polyrigid bundle is polyrigid. That is, it might have deformations as a parabolic bundle but the associated moduli space of parabolic bundles is a point.

\begin{lemma}\label{DIRECTSUM}
Let $\mt$ be a polyrigid\footnote{That is, the corresponding moduli space of parabolic bundles is a point.} semistable bundle (with pairwise difference of weights at any point $<1$). Then, $\mt'=\mt\oplus\dots\oplus\mt$ ($N$-times) is polyrigid for any number of factors.
\end{lemma}

\begin{proof} By induction on length of Jordan-H\"{o}lder series of $\mt$. Let
$\mg=\oplus_{i=1}^k \mv \subset \mt$ be the maximal subbundle with purely factors of the type $\mv$ ($\mv$ is rigid) as in the beginning of this section.

Let $\mg'=\oplus_{a=1}^{N}\mg\subset \mt'.$  

The strategy is somewhat similar to the previous lemma. We show that $\mg'$ deforms with every deformation of $\mt'$ and hence we reduce to  $\tilde{\mt}=\frac{\mt}{\mg}$ (which has generic flags), and we will be done by induction.
The tangent space of $\mg'$ in $\mt'$ is clearly $0$ dimensional so we need
to show that the expected dimension of subbundles of the same Schubert position
as $\mt'$ is $0$. We prove that from knowing that the expected dimension of subbundles with the same Schubert position as $\mg$ in $\mt$ is $0$ (which is evidently so):

Point calculation: If the partial flag variety of $\mt$ at $p_l$ is ($rank(\mt)=n)$, $i^l_1<\dots<i^l_{k_l}=n$ and 
suppose $dim(\mg_{p_l}\cap F_{i^l_a}(\mt_{p_l}))=u^l_a,$
then 
$dim(\mg'_{p_l}\cap F_{Ni^l_a}(\mt_{p_l}))=Nu^l_a,$
so if $Y_l$ is the $l$th Schubert Cell for $\mg$ and $Y'_l$, we have
$codim(Y'_l)=N^2 codim(Y_l)$. So we have from the $0$ expected dimensionality of $\mg$ in $\mt$,

$$rank(\mg)(n-rank(\mg))-degree(\mg)n + degree(\mt)rank(\mg)=\sum codim(Y_l)$$

and the same equation if multiplied by $N^2$, gives the sought for condition on $\mg'$.
\end{proof}
In the proof of the previous theorem we have also shown the following:
\begin{lemma}\label{scale} Let $\mt=\mg_{d,n}$ and suppose $\langle\sigma_{I^1},\dots,\sigma_{I^{s}}\rangle_{d,D,n}\neq 0$, happening in $Gr(r,n)$.
Then we can scale this situation by $N$: Let $J^l$ be subsets of
$\{1,\dots,Nn\}$ each of cardinality $Nr$ given by 
\begin{enumerate}
\item $j^l_{Na}=Ni^l_a$.
\item $j^l_{Na-b}=Ni^l_a-b$ for $0<b<N$.
\end{enumerate}

Then
$\langle\sigma_{J^1},\dots,\sigma_{J^{s}}\rangle_{Nd,ND,nN}\neq 0,$ 
happening in $Gr(Nr,Nn)$.
\end{lemma}

\begin{corollary}\label{PolyRigidInv} Polyrigidity is invariant under Grassmann duality.
\end{corollary}
\begin{proof}
Suppose we are in the situation at the end of section ~\ref{HWitten}. Assume that the Witten bundle of the Grassmann dual $\mq^*$ is polyrigid. Let
$N(\mathcal{M}_2)$ be the Grassmann dual Witten Moduli space  for $f(N)$ (as in section on the saturation conjecture). We know that $N(\mathcal{M}_2)$ is the  set of deformations of $\oplus_{l=1}^{l=N}\mq^*$, hence is a point by polyrigidity (and lemma ~\ref{DIRECTSUM}). Therefore $f(N) =1$ because it is the space of global sections of a bundle over $N(\mathcal{M}_2)$ and hence $\mv$ is polyrigid.
\end{proof}

\begin{lemma}\label{polyrigidinside}
Suppose $\mv\subset\mt$ and corresponds to a polyrigid intersection (where $\mt$ has generic flags at $p_1,\dots,p_s$). Suppose $\mw\subset\mt$ and corresponds to a intersection dimension $0$ situation (that is, there are only finitely many subbundles of $\mt$ with the same Schubert position as $\mw$) and
$\mw\supset \mv$. Notice that the induced flags on $\mw$ are generic. Then,
$\mv\subset\mw$ corresponds to a polyrigid intersection: If the Schubert position of $\mv$ in $\mw$ is $K^1,\dots,K^s$ (with respect to the induced flags), we have
$\langle\sigma_{K^1},\dots,\sigma_{K^s}\rangle_{d,-degree(\mw),rank(\mw)}=1$ and is polyrigid.
\end{lemma}
\begin{proof} Clearly $dim(\mv,\mw)=0$. We need to prove the polyrigidity property. If $\mv\subset\mw$ were not polyrigid, then we will find a $N$ so that (let $rank(\mv)=r_v=r,rank(\mw)=r_w, n=rank(\mt),d_w=-degree(\mw)$)

$$f_1(N)=\langle\sigma_{NK^1},\dots,\sigma_{NK^s}\rangle_{d,N(d_w-d)+d,N(r_w-r)+r}>1,$$

 so we just need to show the `injection' of this number into 

$$f_2(N)=\langle\sigma_{NI^1},\dots,\sigma_{NI^s}\rangle_{d,N(D-d)+d,N(n-r)+r}.$$

The idea is very similar to the arguments so far: let $\mt(N)$=
$\mg_{N(D-d)+d,N(n-r)+r}$, take generic flags on it, and denote the subbundle corresponding to 

$$\langle\sigma_{NI^1},\dots,\sigma_{NI^s}\rangle_{d,N(D-d)+d,N(n-r)+r}=1$$

 by $\mv(N)$. The first thing to check is that the Grassmann  dual Witten bundle to $\mv(N)$ is a deformation of the ``scale by $N$ (see lemma ~\ref{scale})'' of the dual Grassmann Witten bundle of $\mv$ in $\mt$ (together with the parabolic structure). Hence we find a 
$\mw(N)\subset\mt(N)$, $\mw(N)\supset\mv(N)$ corresponding to $\mv_1\supset\mv$. Once again, the induced flags on $\mw(N)$ are generic. This follows from
$dim(\mw(N),\mt(N))-dim(\mw(N)/\mv(N),\mt(N)/\mv(N))=N(dim(\mw,\mt)-dim(\mw/\mv,\mt/\mv))=0$, and $dim(\mw(N)/\mv(N),\mt/\mv(N)=0$, so $dim(\mw(N),\mt(N))=0$,
hence the induced flags on $\mw(N)$ are generic.

It is easy to check that $deg(\mw(N))=-d+N(-(d_w-d))$, and $rank(\mw(N))=r+N(r_w-r)$, so the assumption $f_1(N)>1$, gives us atleast 2 subbundles of $\mw(N)$ of
degree and rank same as that of $\mv$ which have the Schubert position
$NK^1,\dots,NK^s$. We just need to check that as subbundles of $\mt(N)$, they have Schubert positions $NI^1,\dots,NI^s$ so that they contribute to $f_2(N)$ making it larger than 1. This finishes the proof.

\end{proof}

\begin{lemma}\label{lastlemma} Let $\mv\subset \mt$ be from a polyrigid situation of the type:

$$\langle\sigma_{I^1},\dots,\sigma_{I^{s}}\rangle_{d,D,n}=1,$$

Suppose $\mg\subset \mv$ be as produced in the theorem ~\ref{JHH} the Witten Bundle structure on $\mv$. Then, the Witten bundle of $\mg\subset\mt$ is polyrigid (in fact has no deformations as a bundle).
\end{lemma}
\begin{remark} The expected dimension of subbundles of $\mt$ in the Schubert position of $\mg$ is $0$ (this is because $dim(\mg,\mv)=0$, $dim(\mg,\mv)=dim(\mg,\mt)$). 
\end{remark}
\begin{proof}
The proof is really the same as lemma ~\ref{StrongRigid}. It is important at the very outset, to ignore parts of the flag varieties of $\mv$ where the Witten weight does not jump. 
\end{proof}

\subsection{Proof of Theorem ~\ref{UneTheoreme}}\label{TheFinalStep}

Let $(I^1,\dots,I^s,d,D,n)$ be a ``Polyrigid Witten'' situation. Let $\mt=\mg_{D,n}$, $\mv$ the subbundle and $\mq$ the quotient(with the Grassmann dual Witten structure).  Suppose $\mq$ is not semistable. Apply the arguments of the Jordan-H\"{o}lder section to $\mq^*$ and produce a subbundle $\mq_1^*$ with simple factors $\ma$. So we have the following situation (ignoring parabolic structures),

$$0\subset \mv\subset \mv_1\subset \mt$$

with $\mt/\mv=\mq$ and $\mt/\mv_1=\mq_1$.

We also have that 

$dim(\mv_1,\Bbb{C}^n)=0=dim (\mv_1/\mv,\mt/\mv).$

(so $\mv_1$ is a generic bundle too)

From lemma ~\ref{lastlemma}(after dualising) and lemma ~\ref{polyrigidinside}, that $\mv_1$ is polyrigid, and $\mv\subset\mv_1$ corresponds to a polyrigid intersection. Now the Witten bundle of the Grassmann dual of $\mv_1$ may not have only one simple factor. But then we can iterate this procedure. This finishes the proof of theorem ~\ref{UneTheoreme}.

\bibliographystyle{plain}

\def\noopsort#1{}

\end{document}